\documentclass[]{article}

\usepackage[a4paper, margin=2.5cm]{geometry}

\usepackage{authblk}
\usepackage[margin=1cm]{caption}

\usepackage{graphics} %
\usepackage{epsfig} %
\usepackage{amsmath} %
\usepackage{amssymb}  %
\usepackage{booktabs}
\usepackage{makecell}
\usepackage{hyperref}
\usepackage{tikz}

\usepackage{amsfonts}
\usepackage{xcolor}
\usepackage{booktabs}

\usepackage{optidef}

\usepackage{placeins}
\usepackage{comment}
\newcommand{\shortminus}{-}

\usepackage{pythonhighlight}
\lstnewenvironment{mypython}[1][]{\lstset{style=mypython, frame=none, #1}}{}
\usepackage{mathtools}

\definecolor{wheat}{rgb}{0.96,0.87,0.70}
\newcommand{\transp}{^\top}

\newcommand{\ind}[1]{{_{\mathrm{#1}}}}
\newcommand{\diag}[1]{\text{\normalfont diag}\left( {#1} \right)}

\newcommand{\R}{\mathbb{R}}

\newcommand{\fimpl}{f}%

\newcommand{\acados}{\texttt{acados}}

\DeclarePairedDelimiter\abs{\lvert}{\rvert}%
\DeclarePairedDelimiter\norm{\lVert}{\rVert}%

\makeatletter
\let\oldabs\abs
\def\abs{\@ifstar{\oldabs}{\oldabs*}}
\let\oldnorm\norm
\def\norm{\@ifstar{\oldnorm}{\oldnorm*}}
\makeatother

\DeclareMathSymbol{\sm}{\mathbin}{AMSa}{"39}

\newcommand{\dpartial}[2]{\frac{\partial #1}{\partial #2}}

\newcommand{\dtplant}{\Delta t}

\newcommand{\contcontrl}{v}

\newcommand{\nphases}{M}
\newcommand{\ftransk}{\Gamma_k}
\newcommand{\contx}{\xi}
\newcommand{\Tclc}{T\ind{clc}}
\newcommand{\npc}{n\ind{pc}}

\newcommand{\transcontrl}{\eta}

\title{Multi-Phase Optimal Control Problems for Efficient Nonlinear Model Predictive Control with \acados}

\author[1,2]{Jonathan Frey}
\author[1]{Katrin Baumgärtner}
\author[1]{Gianluca Frison}
\author[1,2]{Moritz Diehl}

\affil[1]{Department of Microsystems Engineering (IMTEK), University Freiburg, Germany}
\affil[2]{Department of Mathematics, University Freiburg, Germany}

\date{Preprint. Article submitted to journal \textit{Optimal Control Applications and Methods} on July 12, 2024}

\begin{document}

\abstract{
Computationally efficient nonlinear model predictive control relies on elaborate discrete-time optimal control problem (OCP) formulations trading off accuracy with respect to the continuous-time problem and associated computational burden.
Such formulations, however, are in general not easy to implement within specialized software frameworks tailored to numerical optimal control.
This paper introduces a new multi-phase OCP interface for the open-source software \acados{} allowing to conveniently formulate such problems and generate fast solvers that can be used for nonlinear model predictive control (NMPC).
While multi-phase OCP (MOCP) formulations occur naturally in many applications, this work focuses on MOCP formulations that can be used to efficiently approximate standard continuous-time OCPs in the context of NMPC.
To this end, the paper discusses advanced control parametrizations, such as closed-loop costing and piecewise polynomials with varying degree, as well as partial tightening and formulations that leverage models of different fidelity.
An introductory example is presented to showcase the usability of the new interface.
Finally, three numerical experiments demonstrate that NMPC controllers based on multi-phase formulations can efficiently trade-off computation time and control performance.
}
\thanks{The authors want to thank Rudolf Reiter for fruitful discussions.\\
This research was supported by DFG via Research Unit FOR 2401 and projects 424107692, 504452366, by BMWK 03EN3054B, and by the EU via ELO-X 953348.
}

\maketitle

\section{Introduction}
High-performance algorithms for nonlinear model predictive control (NMPC) %
have extended its real-time applicability from chemical processes to a huge variety of application areas\cite{Dabbene2023, Haenggi2022, Zanelli2021a, Carlos2020, Carlos2020a, Gao2023, Norouzi2022, Romero2022}.
Most NMPC implementations use a discrete-time model with a fixed time step that is equal to the sampling time.
However, more elaborate optimal control problem (OCP) formulations have enormous potential to reduce the computational burden associated with the online optimization. %
Nonuniform time grids, which have been shown to work especially well with accurate cost integration~\cite{Frey2024b}, are an easy first step in this direction.
However, more sophisticated OCP formulations, like piecewise polynomial control parametrizations of varying degree \cite{Vassiliadis1994}, closed-loop-costing formulations \cite{Magni2004}, or using models of different fidelity for parts of the horizon, are rarely used in practice, due to the additional implementation effort.
Such formulations require the problem functions as well as input and state dimensions to vary over the horizon. %
To this end, this paper presents a new feature of the open-source software~\acados{} \cite{Verschueren2021} which allows for a convenient formulation of multi-phase OCPs via its \texttt{Python} interface.

We use the term \textit{multi-phase} OCP (MOCP) to describe OCP formulations that may have structurally different models, constraints and cost formulations on parts of the horizon.
In contrast the term \textit{multi-stage} was used in previous works but is used in different ways, e.g. in the context of tree-structured OCPs \cite{Kouzoupis2019, Fiedler2023, Lucia2014c} and to describe OCP formulations which only allow to vary cost and constraint functions~\cite{MatlabMultiStageNMPC}.
The proposed multi-phase formulation is similar to the one tackled by \texttt{GPOPS-II}~\cite{Patterson2014a}, which is a commercial \texttt{Matlab} software package for solving multi-phase OCPs using state-of-the-art sparse NLP solvers, like \texttt{IPOPT}~\cite{Waechter2006} and \texttt{SNOPT}~\cite{Gill2005}.
\texttt{GPOPS-II} has been successfully used for MOCPs~\cite{Ye2016}.

However, the focus of this paper is on MOCP formulations with OCP-structure exploiting algorithms.
MOCPs have been previously tackled in the direct multiple shooting package MUSCOD-II \cite{Leineweber2003}, which allows multi-phase formulations with piecewise polynomial control parametrizations and coupled processes.
While MUSCOD-II focused on problems from process control with time scales in the second to minute range, \cite{Leineweber2003, Leineweber2003a}, the \texttt{acados} software package has been successfully deployed on systems with much shorter time scales \cite{Haenggi2022, Zanelli2021a, Carlos2020, Carlos2020a, Gao2023, Norouzi2022, Romero2022}.
Both MUSCOD-II and \texttt{acados}~use a direct multiple shooting based formulation.

The open-source software package \acados{} implements efficient algorithms for embedded optimal control \cite{Verschueren2021}.
It is written in \texttt{C} and relies on the linear algebra package \texttt{BLASFEO} which provides performance-optimized routines for small to medium sized matrix operations~\cite{Frison2018}.
The nonlinear OCPs can be solved with \acados~using variants of sequential quadratic programming (SQP), the real-time iteration (RTI)~\cite{Diehl2001} and advanced-step real-time iteration (AS-RTI) scheme~\cite{Frey2024a}, as well as differential dynamic programming (DDP)~\cite{Mayne1966}, which was recently added~\cite{Kiessling2024}.
In \acados, the quadratic subproblems are solved exploiting the block structure of optimal control problems.
It has been shown that OCP structure exploiting solvers can be many times faster compared to general sparse solvers or dense ones used with condensing \cite{Kouzoupis2018, Vanroye2023}.
A variety of quadratic programming~(QP) solvers, such as \texttt{HPIPM}, \texttt{qpOASES}, \texttt{DAQP}, \texttt{OSQP}, \texttt{qpDUNES},
\cite{Frison2020a, Ferreau2014, Arnstrom2022, Stellato2017a, Frasch2013}
are interfaced, which either tackle the OCP-structured QP directly or after applying full or partial condensing to it~\cite{Frison2016, Axehill2015}.
The system dynamics can be handled with integration methods which are able to propagate forward and adjoint sensitivities efficiently in addition to the nominal result~\cite{Verschueren2021, Frey2023}.

The remainder of this introduction aims at giving a brief overview on software frameworks for NMPC other than \acados.
The tool \texttt{CasADi}~\cite{Andersson2019} offers a convenient symbolic framework, algorithmic differentiation and interfaces to many state-the-art NLP solvers, like \texttt{IPOPT}~\cite{Waechter2006} and \texttt{SNOPT}~\cite{Gill2005}.
A variety of software projects specialized on OCP formulations have been developed on top of \texttt{CasADi}s symbolic framework.
The package \texttt{do-mpc}~\cite{Fiedler2023}, provides \texttt{Python} functionality to allow fast prototyping of NMPC, moving horizon estimation (MHE) and supports tree-structured OCP formulations.
In terms of solvers it relies on the ones available in \texttt{CasADi}.
The \texttt{rockit} package~\cite{Gillis2020}, allows to conveniently formulate OCPs in \texttt{Python} and \texttt{Matlab} and also covers MOCP formulations.
In addition to the \texttt{CasADi} solvers, it also provides interfaces to \texttt{acados} and \texttt{fatrop}~\cite{Vanroye2023}, which is an NLP solver highly inspired by \texttt{IPOPT} that exploits the optimal control problem structure.
The \texttt{OpEn}~\cite{Sopasakis2020} software package allows to conveniently formulate single-phase OCP problems and implements the proximal averaged Newton-type method for optimal control (PANOC) as a \texttt{Rust} solver and allows to generate solvers for user defined problems from \texttt{Python} and \texttt{Matlab} using the \texttt{CasADi} symbolics.
The \texttt{GRAMPC}~\cite{Englert2019} software package can generate embedded solvers from \texttt{Matlab} and uses a gradient‑based augmented Lagrangian method.
The commercial solver \texttt{FORCESPRO}~\cite{Zanelli2017b} implements competitive algorithms for NMPC which support varying dimensions between the stages~\cite{Zanelli2017b, Domahidi2012}.
However, benchmark results created with the academic license of \texttt{FORCESPRO} could not be disclosed in past research~\cite{Kouzoupis2018, Norouzi2022}.

\textbf{Outline.}
The paper introduces the multi-phase OCP formulation in Section~\ref{sec:ocp_formulations} and discusses how it can be handled using direct multiple shooting.
Section~\ref{sec:multiphase_examples} motivates using advanced OCP formulations which may be cast as MOCPs to design efficient NMPC controllers.
In particular, it discusses piecewise polynomial control parameterizations, partial and progressive tightening formulations, and the use of models of different fidelity within a single OCP.
Moreover, it presents and conceptually extends closed-loop costing formulations.
Section~\ref{sec:implementation} discusses the efficient treatment of multi-phase formulations within the \acados~software package and presents a tutorial example.
Section~\ref{sec:experiments} demonstrates how NMPC controllers based on the MOCP formulations detailed in Section~\ref{sec:multiphase_examples} are able to trade off computation time and control performance in ways that are not accessible when one is limited to single-phase OCPs.

\section{Multi-Phase Optimal Control Problem Formulations}
\label{sec:ocp_formulations}
This section presents a continuous-time multi-phase optimal control problem (MOCP) formulation in Section~\ref{sec:multiphase_formulation} and describes how it can be discretized using direct multiple shooting in Section~\ref{sec:mocp_shooting}.

\subsection{Continuous-time multi-phase optimal control problem}
\label{sec:multiphase_formulation}
Throughout this paper, we will treat continuous-time multi-phase optimal control problems (MOCP), which can be stated as
\begin{mini!}
	{\substack{\contx_1(\cdot), \dots, \contx_{\nphases}(\cdot), \\
			\contcontrl_1(\cdot), \dots, \contcontrl_{\nphases}(\cdot), \\
			\transcontrl_1, \dots, \transcontrl_\nphases}}
	{\sum_{k=1}^{\nphases}\int_{t_k}^{t_{k+1}}\ell_k(\contx_k(t), \contcontrl_k(t)) \, \mathrm{d} t + E_k(\contx_k(t_{k+1}), \transcontrl_k)
		\label{eq:MOCP_cost}}
	{\label{eq:MOCP}}
	{}
	\addConstraint{0}{ = \bar{x}_0 - \contx_{1}(t_1) \label{eq:MOCP_initial_state}}
	\addConstraint{0}{= \fimpl_k(t, \contx_k(t), \dot{\contx}_k(t), \contcontrl_k(t))\label{eq:MOCP_ODE}}
	\addConstraint{0}{\geq g_k(\contx_k(t), \contcontrl_k(t))}{}
	\addConstraint{0}{=\ftransk(\contx_k(t_{k+1}), \transcontrl_k) - \contx_{k+1}(t_{k+1}) \label{eq:MOCP_transition}}
	\addConstraint{0}{\geq g\ind{e}(\contx_{\nphases}(T))}
	\addConstraint{\quad \mathrm{for}\; t\in[t_k, t_{k+1}),\; k=1,\dots, \nphases. \nonumber }
\end{mini!}
The finite time horizon $[0, T]$ is split into $\nphases$ fixed subintervals $[t_k, t_{k+1}]$ with $t_1 = 0$ and $t_{M+1} = T$.
Each interval defines a phase $k\in\! \{1, \dots, \nphases \}$.
For each phase $k$, an implicit ODE $f_k$
is given, which defines the state trajectory $\contx_k: [t_k, t_{k+1}]\rightarrow \R^{n_{x,k}}$ for a given control trajectory $\contcontrl_k: [t_k, t_{k+1}]\rightarrow \R^{n_{\contcontrl,k}}$ and initial state.
The initial state of the first phase is given by $\bar{x}_0$, while the initial state for the subsequent phases is given by the transition functions $\ftransk: \R^{n_{x, k}} \times \R^{n_{\transcontrl, k}} \rightarrow \R^{n_{x, k+1}}$, which map the terminal state of phase $k$ and discrete decision variables $\transcontrl_k \in \R^{n_{\transcontrl,k}}$ to the initial state of the next phase, allowing for the state dimension to vary between the phases.
The functions $\ell_k$ and $E_k$ define the path and terminal cost terms for phase $k$.
The functions $g_k$ summarize the inequalities imposed on states and controls in phase $k$.
Additionally, the function $g\ind{e}$ summarizes the terminal constraints imposed at the end of the horizon.

The transition formulation \eqref{eq:MOCP_transition} is similar to the one in MUSCOD-II \cite{Leineweber2003b} with the addition of the discrete decision variable~$\transcontrl_k$.
Additionally, the formulation in the work by Leineweber et al. \cite{Leineweber2003b} considers global variables and the time grid points $t_i$ as optimization variables.
The global variables are implemented as separate variables on each shooting interval and are constrained to be equal. %
Global variables can be formulated in \eqref{eq:MOCP} by state augmentation on the full horizon.
In order to have time grid points as optimization variables, the dynamics can be augmented with a clock state and a speed of time variable acting as a control.
The differences in problem formulation can be motivated by the different solution methods implemented in \acados~and MUSCOD-II, respectively.
While {MUSCOD-II} deploys sparse linear algebra on blocks and allows linear couplings between stages, \texttt{acados} uses specialized algorithms for purely OCP-structured problems in which the dynamics constraints are the only coupling between stages.

MOCP formulations occur naturally in different application when formulating OCPs where the dynamic behavior qualitatively changes at a certain point in time.
For example, when considering walking robots \cite{Schultz2009},
chemical plants, where some amount of substance is added at a certain point in time, e.g. when considering recycled waste cuts~\cite{Diehl2002e},
or multi-train scheduling~\cite{Ye2016}.
However, the main focus of this paper is on MOCP formulations which can be derived to approximate the solution of a continuous-time infinite-horizon problem in the context of NMPC as discussed next.

\subsection{Multi-phase multiple shooting discretization}
\label{sec:mocp_shooting}
One can discretize MOCP \eqref{eq:MOCP} using direct multiple shooting \cite{Bock1984} dividing each time interval $[t_k, t_{k+1}]$ into $N_k$ shooting intervals.
The resulting nonlinear program (NLP) can be written as
\begin{mini!}
    {\substack{\boldsymbol{x}, \boldsymbol{u}, \boldsymbol{\transcontrl}}}
    {\sum_{k=1}^{\nphases }  \sum_{j=0}^{N_k\shortminus 1}   L_{k,j}(x_{k,j}, u_{k,j}, \Delta t_{k, j})  +  E_k(x_{k, N_k}, \transcontrl_k)
    \label{eq:dmocp_cost}}
    {\label{eq:dmocp}}
    {}
    \addConstraint{x_{1,0}}{= \bar{x}_0 \label{eq:dmocp_initial_state}}
    \addConstraint{x_{k+1, 0}}{=\ftransk(x_{k,N_k}, \transcontrl_k) %
    \label{eq:dmocp_transition},}{\,}{k = 1,\dots M\shortminus 1}
    \addConstraint{x_{k, j+1}}{= \phi_{k,j}(x_{k,j}, u_{k,j}),\,\;} {j=0,\dots, N_k \shortminus 1 \label{eq:dmocp_dynamics}}{, k = 1,\dots M}
    \addConstraint{0}{\geq g_{k,j}(x_{k,j}, u_{k,j}),} {j=0,\dots, N_k\shortminus 1}{, k = 1,\dots M}
    \addConstraint{0}{\geq g\ind{e}(x_{\nphases, N_\nphases}),}
\end{mini!}
where $ \boldsymbol{x}=(x_{1,0},\dots, x_{1, N_1}, x_{2, 0}, \dots, x_{\nphases, N_\nphases})$, $\boldsymbol{u}=(
u_{1,0}, \dots, u_{1, N_1-1}, u_{2, 0}, \dots, u_{\nphases, N_\nphases \shortminus 1})$, and $\boldsymbol{\transcontrl}=(\transcontrl_1, \dots, \transcontrl_\nphases)$.
The discrete values $x_{k,j} \in \R^{n_{x,k}}$ represent the values of $\contx_k(\cdot)$ at the corresponding shooting nodes.
The controls $u_{k,j} \in \R^{n_{u,k}}$ act on the $j$th shooting interval of phase $k$.
The cost functions $L_{k,j}$ are called stage costs and reflect the cost integrated over the shooting interval $[t_{k,j}, t_{k, j+1})$ with $\Delta t_{k, j} = t_{k, j+1} - t_{k, j}$.
The initial state of the first phase is given by \eqref{eq:dmocp_initial_state} and for subsequent stages by \eqref{eq:dmocp_transition}.
The discrete-time dynamics are described by the functions $\phi_{k,j}$, which represent an integration scheme applied to the continuous-time dynamic system in \eqref{eq:MOCP_ODE}.
Dedicated functions that compute the numerical integration of continuous-time dynamics and the solution sensitivities are an essential component of for efficient solution of OCPs and available in \acados~\cite{Frey2023}.

Finally, the constraint functions $ g_{k,j}$ represent the constraints $g_k$ on shooting interval $[t_{k,j}, t_{k, j+1})$.
Most commonly, the constraints are only enforced at the initial point of the shooting interval, i.e. for $x_{k,j}, u_{k,j}$.
However, it is possible to also impose them on intermediate points.
Similarly, the stage cost $L_{k,j}$ often corresponds to a simple Euler integration of the continuous-time cost $\ell_k$ over a shooting interval.
Especially, when using longer intervals, higher order integration of the cost term are necessary to retain a good approximation quality.
The cost integration can be performed efficiently together with the integration of the dynamics~\eqref{eq:MOCP_ODE}.
More details on this can be found in Section~\ref{sec:cost_integration}.

\section{Multi-phase OCP formulations for Model Predictive Control}
\label{sec:multiphase_examples}
This section motivates the use of multi-phase optimal control problem (MOCP) formulations for NMPC applications, in which an ideal controller would apply the exact solution to a continuous-time infinite-horizon optimal control problem, which is presented in Section~\ref{sec:motivation}.
Moreover, Section~\ref{sec:motivation} introduces the concept of the cost-to-go function, the importance of its approximation in NMPC, and why MOCP formulations should be considered in this regard.

The following subsections provide examples of optimal control problem formulations which are suitable to approximate parts of the cost-to-go term corresponding to different parts of the time horizon in successively more approximate ways.
These formulations can be cast as multi-phase problems and can thus be solved efficiently using \acados.
In particular, Section~\ref{sec:models_of_different_fidelity} motivates using models of different fideltiy for different parts of the horizon, Section~\ref{sec:pw_polynomial} discusses piecewise polynomial control parameterizations, Section~\ref{sec:clc} discusses closed-loop costing and Section~\ref{sec:ptight} presents partial tightening.

\subsection{Continuous-time optimal control problem, model predictive control and cost-to-go} \label{sec:motivation}
The standard infinite-horizon continuous-time optimal control problem (OCP) which we aim at approximating in various ways in this paper can be written as
\begin{mini!}|s|
	{\substack{x(\cdot), \contcontrl(\cdot)}}
	{\int_{0}^{\infty} \ell(x(t), \contcontrl(t)) \, \mathrm{d} t
		\label{eq:NOCP_cost}}
	{\label{eq:NOCP}}
	{{V(\bar{x}_0)} = }
	\addConstraint{\bar{x}_0}{= x(0) \label{eq:NOCP_initial_state}}
	\addConstraint{0}{= \fimpl(x(t), \dot{x}(t), \contcontrl(t)),~\;}{t\!\in\![0, \infty) \label{eq:NOCP_ODE}}
	\addConstraint{0}{\geq g(x(t), \contcontrl(t)),}{t\!\in\![0, \infty),}
\end{mini!}
where $ x: [0, \infty) \rightarrow \R^{n_x}$, $\contcontrl: [0, \infty) \rightarrow \R^{n_\contcontrl} $
are the state and control trajectories respectively, $ \bar{x}_0 $ is the initial state value, the function $ \fimpl $ describes the implicit system dynamics and $ g $ denotes the inequality constraints.
The optimal value of OCP \eqref{eq:NOCP} is defined as $V(\bar{x}_0)$.
The function $V(\cdot)$ is called cost-to-go and corresponds to the value function in the field of reinforcement learning.

The goal of model predictive control (MPC) is to operate a system, which typically applies a constant control input for a fixed sampling time $\dtplant$.
This practical constraint can be formalized as
\begin{align}
	\contcontrl(t) = u_0 \quad \mathrm{for} \; t \in [0, \dtplant]. \label{eq:pw_const_plant}
\end{align}

Regarding OCP~\eqref{eq:NOCP} from a dynamic programming point of view allows us to split the infinite horizon in different parts which might be approximated in different ways.
Using the principle of optimality and accounting for \eqref{eq:pw_const_plant}, we can rewrite \eqref{eq:NOCP} as
\begin{mini}
	{\substack{x(\cdot),u_0} }
	{\int_{0}^{\dtplant} \ell(x(t), u_0) \, \mathrm{d} t + V(x(\Delta t)) \label{eq:DP_OCP}}{}{}
	\addConstraint{\bar{x}_0}{= x(0)}
	\addConstraint{0}{= \fimpl(x(t), \dot{x}(t), u_0),~}{t\! \in [0, \dtplant]}
	\addConstraint{0}{\geq g(x(t), v(t)),}{t\! \in [0, \dtplant)}
	,
\end{mini}
which is again an OCP but with a much shorter horizon of length $\Delta t$.
Of course, the whole complexity of the problem is shifted into the minimization of the cost-to-go term with this reformulation.
However, this shows that approximating the cost-to-go is key in the developement of an efficient NMPC controller.

In order to obtain an OCP formulation for MPC, one typically selects a finite time horizon, replaces the infinite integral in \eqref{eq:NOCP} with a finite one and adds a terminal cost on the terminal state.
This approximation of the cost-to-go together with \eqref{eq:DP_OCP} can be discretized using multiple shooting, where the first shooting interval typically corresponds to~\eqref{eq:DP_OCP}.

Additionally, the control trajectory to approximate the closed-loop cost could be parameterized in different ways.
In particular, piecewise polynomial control parameterizations and closed-loop costing are dicussed in Section~\ref{sec:pw_polynomial} and Section~\ref{sec:clc}.
These parameterizations correspond to piecewise continuous functions, that are not applicable to the real system, given the practical constraint \eqref{eq:pw_const_plant}.
However, since the cost-to-go is anyway approximated, such control parametrizations can lead to better approximations of the cost-to-go when using longer intervals.

In the following, we derive and discuss various multiple shooting based approximate versions of \eqref{eq:NOCP} which apply structurally different approximations for parts of the infinite horizon in~\eqref{eq:NOCP} and can be be phrased as MOCPs.

\subsection{Models of different fidelity}\label{sec:models_of_different_fidelity}
In many physical systems, models of different fidelity are available and choosing an appropriate one for NMPC might depend on the control frequency, desired time horizon and the available solvers.
The consideration of high-fidelity models is typically beneficial for MPC performance.
However, such models might be computationally demanding and can become unstable for long time horizons, which are important to consider in OCPs in order to capture the evolution of slow dynamics and economic cost terms.
Multi-phase OCP formulations allow one to choose different models for specific parts of the horizon, allowing for many more degrees of freedom when deriving a discrete-time finite horizon problem.
Low-fidelty models are typically cheaper to integrate with a certain accuracy, not only due to their potentially reduced dimensionality, but also because they are usually less stiff and thus can be handled by computationally cheaper integration schemes, such as explicit Runge-Kutta (ERK) methods.
On the other hand, high fidelty models typically comprise both fast and slow modes, which causes them to be stiff.
Such stiff models need to be handled with computationally expensive implicit Runge-Kutta (IRK) methods or a great number of very small integration steps of ERK methods.

Moreover, some nonlinear constraints, which are expensive to evaluate, could be only included in a first part of the horizon.
This can be motivated by the practical observation that constraints tend to be active in the first part of the horizon \cite{Frey2020}.
For example, modelling oscillations within a physical system can be expensive and only meaningful on short time horizons, as in the application of wind turbine control~\cite{Schlipf2013}.

Lastly, a variety of NMPC applications rely on underlying controllers that handle the actuators.
However, it can be beneficial to directly model these actuators to accurately represent their behavior at least in the first part of the horizon.
This can allow one to specify cost and constraints that need a model of the underlying actuator.
Including such considerations can unleash optimality potential that is inaccessible otherwise.

\subsection{Piecewise polynomial control parameterization}
\label{sec:pw_polynomial}
While piecewise constant control discretizations are by far the most common parameterization for multiple shooting based MPC, formulation \eqref{eq:dmocp} can accommodate piecewise polynomial control parameterizations with degrees varying between shooting intervals.
The $k$-th component of the control trajectory on a shooting interval $[\tau_0, \tau\ind{e}]$ could be parameterized by a polynomial
\begin{align}
	v_{n,k}(t) = \sum_{i=0}^{n\ind{deg}} u_{n,k,i} \cdot (t-{\tau_0})^i
\end{align}
of degree $n\ind{deg}$, such that the discrete control input is given by 
$u_n = (u_{n,k,i})_{{i=0,\dots, n\ind{deg}, \,k=1,\dots, n_\contcontrl}} $.
This parametrization offers more degrees of freedom and can in general better approximate the optimal continuous-time control trajectory. %
Higher order control parametrizations could thus allow the use of longer shooting intervals compared to a piecewise constant parametrization, while maintaining the approximation quality of the MPC feedback law. %

The practical MPC consideration in \eqref{eq:pw_const_plant} motivates using a constant control parameterization on the first interval $[0, \dtplant]$.
In order to combine a constant control input on the first shooting interval with higher order control parameterizations on other shooting intervals, an MOCP formulation is thus necessary.

For linear parametrizations, simple bounds on the inputs  can be satisfied everywhere by enforcing them at the start and end of each control interval.
In contrast, for higher order polynomial parametrizations, even simple control bounds might be violated within the shooting intervals if only enforced at the boundary points.
In order to avoid excessive use of such violations, one possibility is to enforce the control bounds at $ \npc $ equidistant points within every shooting interval.
For any fixed point $\bar{\tau} \in [t_n, t_{n+1}] $, we have
\begin{align}
v_{n,k}(\bar{\tau}) = \sum_{i=0}^{n\ind{deg}} u_{n,k,i} \cdot (\bar{\tau}-{t_n})^i.
\end{align}
which results in a linear inequality constraint with coefficients $(\bar{\tau}-{t_n})^i$ for every intermediate point on which a control bound is enforced.
Additionally, one could add smooth penalties on violations of the control bounds and integrate them over the shooting intervals, see Section~\ref{sec:cost_integration}.

\subsection{Cost-to-go approximation via closed-loop costing}
\label{sec:clc}
The idea of \textit{closed-loop costing} (CLC) is to apply a simple control law $\kappa: \R^{n_x} \rightarrow \R^{n_\contcontrl}$ to the nonlinear system in \eqref{eq:NOCP_ODE} and assign the resulting cost term as a terminal cost\cite{Nicolao1998, DeNicolao1996a, Diehl2004f, Magni2004}.
This control law is typically a locally stabilizing linear controller, e.g. corresponding to an LQR controller.
In the Reinforcement Learning and Dynamic Programming community, this idea is more commonly referred to as rollout of a base policy\cite{Bertsekas1996a}.
The cost-to-go is approximated by the integrated cost associated with the closed-loop system under policy $\kappa$.

The closed-loop costing problem corresponding to \eqref{eq:NOCP} can be written as
\begin{mini!}|s|
    {\substack{x(\cdot), \contcontrl(\cdot)}}
    {\int_{0}^{T_1} \ell(x(t), \contcontrl(t)) \, \mathrm{d} t + \int_{T_1}^{\Tclc} \ell(x(t), \kappa(x(t))) \, \mathrm{d} t %
    \label{eq:approx_clc_ocp_cost}}
    {\label{eq:approx_clc_ocp}}
    {}
    \addConstraint{0}{= x(0)- \bar{x}_0 \label{eq:approx_clc_ocp_initial_state}}
    \addConstraint{0}{= \fimpl(t, x(t), \dot{x}(t), \contcontrl(t)),}{t\!\in\![0, T_1) \label{eq:approx_clc_ocp_ODE_1}}
    \addConstraint{0}{= \fimpl(t, x(t), \dot{x}(t), \kappa(x(t))),~}{t\!\in\![T_1, \Tclc) \label{eq:approx_clc_ocp_ODE_2}}
    \addConstraint{0}{\geq g(x(t), \contcontrl(t)),}{t\!\in\![0, T_1]}
    \addConstraint{0}{\geq g(x(t), \kappa(x(t))),}{t\!\in\![T_1, \Tclc)}
    \addConstraint{0}{= g\ind{e}(x(\Tclc))},
\end{mini!}
where the first part of the horizon $[0, T_1]$ is the so-called \textit{control horizon} and the latter part the \textit{simulation horizon} or \textit{closed-loop costing horizon}.
Note that the control input is only defined on the control horizon, i.e. $v:[0, T_1] \rightarrow \R^{n_v}$ and the control dimension is zero on the CLC horizon.
A discrete-time result showing that the stability region grows when increasing the CLC horizon is presented in the work by Magni and colleagues~\cite{Magni2001}.
Note that the constraints are still imposed on the CLC horizon and violations correspond to infinite cost values \cite{Nicolao1998}.

Closed-loop costing problems can be expressed using two phases.
In the terminal phase, the controls are replaced by the control law in the cost and dynamics expressions, resulting in a phase $k$ with $n_{u,k} = 0$.
In literature, single shooting has been the method of choice to handle the CLC horizon due to its superior computational efficiency in the case $n_u=0$.
However, the beneficial convergence properties also make multiple shooting an attractive option in \acados: the partial condensing algorithm provided by \texttt{HPIPM} allows to condense blocks of an arbitrary number of stages, and this can be exploited to condense away all shooting nodes in the CLC horizon, while possibly retaining them in the control horizon.
This makes the computational efficiency of multiple shooting similar to the one of single shooting, as the underlying QP solver does not see the shooting intervals corresponding to the CLC horizon.
For the experiments in Section~\ref{sec:experiments_pendulum_mocp}, it was necessary to use CLC with multiple shooting to achieve convergence.
To the best of our knowledge, the approximate infinite horizon corresponding to the closed-loop costing phase was only implemented with a single shooting interval in previous works \cite{Nicolao1998, DeNicolao1996a, Diehl2004f, Magni2004, Quirynen2014a}.

\subsubsection{Practical treatment of constraints with closed-loop costing}
In practical OCP formulations state constraints are often replaced by penalties in the cost function to avoid infeasible OCPs~\cite{Rawlings2017}.
Especially, when using longer shooting intervals, it is very important to accuaretely integrate the cost term, including penalties~\cite{Frey2024b}.

Accurate cost integration together with a propagation of the cost gradient and a Gauss-Newton Hessian can be performed efficiently.
This is realized by the Gauss-Newton Runge-Kutta (GNRK) integrators and implemented in \acados~\cite{Frey2024b}.
They can effectively handle $L_2$ constraint penalties.
\label{sec:cost_integration}
This implementation has been extended to treat more general convex-over-nonlinear cost terms of the form
\begin{align}
	\ell(x(t), \contcontrl(t)) = \phi(r(x(t), \contcontrl(t)))
\end{align}
with a smooth convex function $\phi$ and a nonlinear function $r $.
The new implementation is able to integrate such cost terms together with their generalized Gauss-Newton (GGN) Hessian \cite{Verschueren2016, Baumgaertner2022}.
This is required to handle more general smooth penalties effectively, such as the squashed barrier closed-loop costing formulation described in Section~\ref{sec:squahed_clc}.

Despite the reformulation of state constraints as penalties, even simple control bounds can render OCP \eqref{eq:approx_clc_ocp} infeasible if the policy $\kappa(\cdot)$ would choose controls that violate those bounds on the CLC horizon.
This would implicitly define a terminal region constraint.
If one wants to avoid that, a blunt, but practical approach would be to not impose the control bounds on the CLC horizon, when using a linear control law $\kappa$.
Another option is to replace control bounds on the CLC horizon with additional penalties.
A third option is to combine the penalty approach with "squashing", as detailed next.

\subsubsection{Squashed closed-loop costing}
\label{sec:squahed_clc}
The idea of squashed CLC is to create a CLC OCP which is aware of the control bounds on the CLC horizon.
This can be achieved by defining the CLC control law as a squashed version of the linear control law, ensuring that control bounds are always satisfied.
A function $\sigma: \R \rightarrow (-1, 1) $ is called \textit{squashing function}, if it is twice continuously differentiable, strictly monotonically increasing, odd, and satisfies $ \lim_{z\rightarrow \pm \infty} \sigma(z) = \pm 1 $.
Additionally, it is important, that $\dpartial{\sigma}{z}(0) = 1$, such that a linear control law is locally not changed.
One example is $\sigma(z) = \tanh(z) $, which was also used in previous works\cite{Baumgaertner2022a}.

In particular, a control variable $u_i$ with symmetric bounds $[-\bar{u}_i, \bar{u}_i] $, can be reformulated using a new variable $v$ and replacing $u_i$ by the expression $ \bar{u}_i \cdot \tanh(\frac{v}{\bar{u}_i}) $ in the OCP.

It is recommended to use squashing together with barrier penalties such that an iterative solver does not go arbitrarily close to the squashed boundaries and gets stuck there\cite{Baumgaertner2022a}.
A function $\beta: \R \rightarrow \R_{\geq 0} \cup \{ \infty \} $ is called barrier function if it is twice continuously differentiable, strictly convex and $\lim_{z\rightarrow \pm 1} \beta(z) = \infty $.
In particular, $\beta(z) = -\log(1+z) - \log(-z+1) $ is such a function and is used in this paper.
The squashed LQR CLC variant with a progressively increasing barrier, can be interpreted in the framework of progressive tightening \cite{Baumgaertner2023b}, which offers both, numerical benefits and closed loop stability guarantees for constrained nonlinear MPC.
It is also possible to extend this concept to one-sided constraints by using a smooth-max function instead of the sigmoid as a squashing function and a one sided barrier.

\subsection{Partial tightening} \label{sec:ptight}
The partial tightening concept was introduced by Zanelli et al.\cite{Zanelli2017a} as a strategy to reduce the computational cost of solving NMPC problems.
To this end, an approximate formulation is introduced which divides the horizon into two phases.
On the second phase, the tightened horizon, the constraints are tightened, i.e. replaced by barrier formulations.
More precisely, a constraint $0\geq g_j(x(t), v(t))$ as in~\eqref{eq:NOCP} is replaced by a logarithmic barrier term
\begin{align}
-\tau \log(g_j(x(t), v(t))), \label{eq:barrier}
\end{align}
which is added to the cost function with a barrier parameter $\tau$.

The partial tightening formulation allows to derive closed-loop stability guarentees, based on the property that the stage cost are monotonically increasing with the stage index, \cite{Zanelli2017a}.
This concept has been further extended to \textit{progressive tightening}, which extends the concept to e.g. formulations with stage-wise increasing barrier parameters $\tau$.

Partial and progressive tightening formulations have been analyzed in the literature regarding closed-loop stability properties \cite{Zanelli2021c, Baumgaertner2023b} and successfully applied in practical applications such as collision avoidance for motion planning \cite{Reiter2024a}, control of a human-sized quadrotor \cite{Zanelli2018}, as well as voltage control in active distribution networks \cite{Valverde2013}.

The barrier formulation in \eqref{eq:barrier} can be combined with a series of monotonically increasing values for $\tau$ for subsequent stage costs on the tightened horizon to a progressive tightening formulation.

In context of the Real-Time Iteration (RTI), the variables corresponding to this second phase can be fully eliminated in the preparation phase using a Riccati recursion.
In an SQP setting this can be achieved by performing partial condensing and summarizing all stages of the tightened horizon into a single block of the reduced QP.

\section{Efficient implementation in \acados}
\label{sec:implementation}
This section presents how MOCPs can be posed for an efficient solution in \acados.
Starting from problem~\eqref{eq:dmocp}, we discuss in Sections~\ref{sec:transition} and~\ref{sec:mocp_acados} the treatment of transition stages and how the problem fits into the internal formulation treated by \acados.
Section~\ref{sec:tutorial_example} presents a tutorial MOCP example and how it can be formulated for an efficient solution in \acados.

\subsection{Treating transition stages}\label{sec:transition}
In order to incorporate transitions between two phases into an existing SQP software for OCP structured problems, we first consider the simple case where \eqref{eq:dmocp_transition} boils down to $ x_{k+1, 0} = x_{k, N_k} $, $n_{\transcontrl,k} = 0$ and $E_k(\cdot) \equiv 0$.
These trivial phase transitions can be eliminated by enforcing
$\phi_{k, N_k\shortminus 1}(x_{k, N_k\shortminus 1}, u_{k, N_k\shortminus 1} ) = x_{k+1, 0} $ and removing $x_{k, N_k}$, $\transcontrl_k$ from the problem.

In the nontrivial case, we observe that the transition equation~\eqref{eq:dmocp_transition} has the same form as the equation for a shooting gap~\eqref{eq:dmocp_dynamics}.
We thus call these \textit{discrete transition stages}.
On the discrete transition stage, the discrete decision variable $\transcontrl_k$ takes the role of a control input and $E_k$ takes the role of the stage cost.

As opposed to formulating the transition as a concatenation of a dynamics step and a transition, the implementation of a transition with a separate discrete transition stage has the following advantages:
Both the terminal state of the phase before the transition and the initial state of the phase after the transition are readily available in the solver and the transition cost $E_k$ fits seamlessly into an MOCP formulation.

\subsection{Multi-phase multiple shooting formulation in \acados}
\label{sec:mocp_acados}

The discrete-time MOCP in~\eqref{eq:dmocp} can be framed as a regular OCP-structured problem with stage-varying costs and constraints and where -- in contrast to the standard setting -- the control and state dimension might vary stage-wise,
\begin{mini!}
	{\substack{x_0,\dots, x_{N},\\
			u_{0}, \dots, u_{N\shortminus 1}}}
	{\sum_{j=0}^{N}L_{j}(x_j, u_j, \Delta t_j) + E(x_N)
	}
	{\label{eq:ocp}}
	{}
	\addConstraint{\bar{x}_0}{= x(0) }
	\addConstraint{x_{j+1}}{= \phi_j(x_j, u_j),~j=0,\dots, N\shortminus 1 \label{eq:acados_mocp_eq}}
	\addConstraint{0}{\geq h_j(x_j, u_j),~j=0,\dots, N\shortminus 1.}
\end{mini!}
This OCP consists of $N$ stages, which capture both the shooting intervals and the transition stages.
In \acados, nontrivial transitions are modeled by adding an extra phase comprised of a single interval.
The transition function $\ftransk $ is specified using the \texttt{acados} discrete dynamics module and the discrete decision variable $\transcontrl_k$ corresponds to the control variable for this phase.
The number of stages $N$ for the OCP~\eqref{eq:ocp} would thus be $N = M-1 +\sum_{k=1}^M N_k$ if all transitions are nontrivial and $N = \sum_{k=1}^M N_k$ if all transitions are trivial.
Note that the stage cost $L_j$ does not depend on the shooting interval length if $j$ is a transition stage.

A major challenge in tackling OCP~\eqref{eq:ocp} with efficient structure-exploiting solvers is that the dimensions of states, controls and constraints are varying arbitrarily between stages.
In \texttt{acados}, this can be easily achieved by the different modules of the SQP type algorithm.
Each stage has its own module to compute and linearize cost and constraint functions.
Each gap constraint~\eqref{eq:acados_mocp_eq} is associated with a module to evaluate and linearize it.
The \texttt{acados} internal OCP-structured QP subproblem is based on the \texttt{HPIPM} software package~\cite{Frison2020a}.
In addition to interior-point solvers for dense and OCP-structured QP formulations, \texttt{HPIPM} offers efficient routines for transforming OCP-structured QPs into dense QPs via full condensing or smaller OCP-structured QPs via partial condensing~\cite{Frison2015a}.
The \texttt{HPIPM} core algorithms, full and partial condensing, and the Riccati-based QP solver, all support stage-varying dimensions.
These full and partial condensing algorithms allow any dense QP solver, such as \texttt{DAQP} \cite{Arnstrom2022} and \texttt{qpOASES} \cite{Ferreau2014} or \texttt{HPIPM} itself, as well as the OCP-structure exploiting solver in \texttt{HPIPM}, to be used seamlessly for multi-phase problem formulations in \acados.

\subsection{Tutorial example of a multi-phase OCP}
\label{sec:tutorial_example}

\begin{figure}[]
	\centering
	\hfill
	\includegraphics[height=6cm]{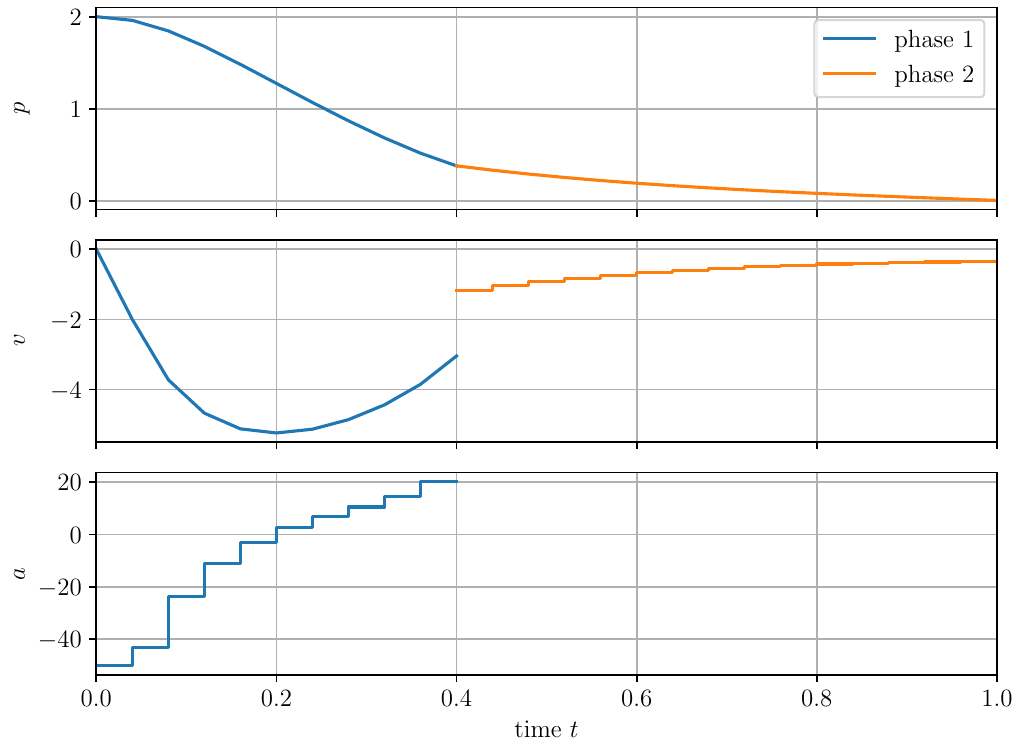}
	\hfill
	\includegraphics[height=6cm]{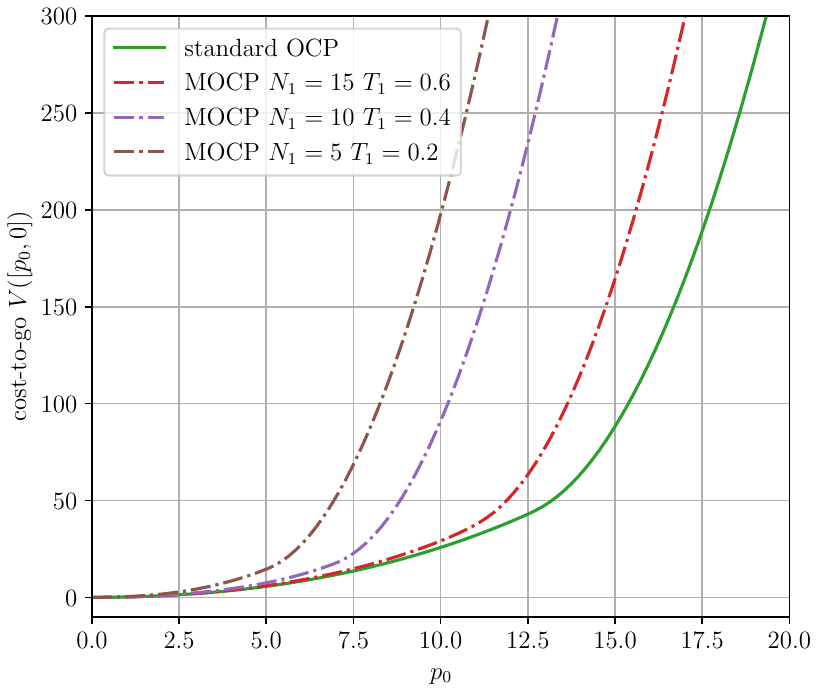}
	\hfill
	\caption{Tutorial example: We consider a double integrator, which is in the second phase approximated by a single integrator model.
		The left plot shows the state and control trajectories for the tutorial MOCP.
		The right plot shows the cost-to-go function as a function of the initial position $p_0$ for different approximations of the OCP which only uses the double integrator model.
		\label{fig:tutorial_transition}
	}
\end{figure}

In this section, we discuss a tutorial example of an MOCP which uses models on distinct parts of the horizon and detail how to formulate the problem in \acados.
The code to reproduce the figures in this section is publicly available\footnote{\url{https://github.com/acados/acados/blob/master/examples/acados_python/mocp_transition_example/main.py}}.
We regard a double integrator system which consists of the state $\contx_1 = (p, s) $, with position $p$ and speed~$s$.
The control input $\contcontrl_1$ is the acceleration $a$ and the continuous-time dynamics are simply given by

\begin{align*}
	0
	&=
	\fimpl_1(\contx_1, \dot{\contx}_1, \contcontrl_1) =
	\begin{bmatrix}
		\dot{p} - s \\ \dot{s} - a
	\end{bmatrix}\!.
\end{align*}
Let us formulate an MOCP with two phases where an approximate model is employed in the second phase.
The approximate model does not consider acceleration and regards the velocity as the control input instead yielding a one dimensional system, i.e.,
$\contx_2 = p$, $\contcontrl_2 = s$.
The system dynamics are then
$\fimpl_2(\contx_2, \dot{\contx}_2, \contcontrl_2) = \dot{p} - s$.
The transition function is $\ftransk(\contx_1) = p$.
We define an MOCP with a time horizon of $T=1$ and the cost functions
\begin{align*}
	\ell_1(\xi_1, \contcontrl_1) &= p^2 + 0.1 s^2 + 10^{\shortminus 3} a^2, &
	E_1(\xi_1) &= p^2 + 0.01 s^2,  \\
	\ell_2(\xi_1, \contcontrl_2) &= p^2 + 0.1 s^2, &
	E_2(\xi_2) &= 10 p^2.
\end{align*}
We impose control bounds as constraints on both phases:
In the first phase, we have $ -50 \leq a \leq 50$, in the second phase, we use $ -5 \leq s \leq 5 $.
We associate the first phase with the interval $ [0, 0.4] $ and the second phase with $ [0.4, 1.0] $ and divide the both intervals uniformly into 10 and 15 shooting intervals, respectively.
The resulting problem is a QP and solved with a single SQP iteration in \acados.
The optimal trajectory for $\bar{x}_0 = [2, 0]\transp$ is visualized in the left plot in Figure~\ref{fig:tutorial_transition}.

The right plot in Figure~\ref{fig:tutorial_transition} shows how MOCPs can approximate the cost-to-go function.
All OCPs use $N=25$ shooting intervals and a time horizon of $T=1$.
The MOCPs use the single integrator approximate OCP phase as described above on the interval $[T_1, T]$.
The shooting intervals are equidistant on both phases with $N_1$ intervals in the first phase and $N-N_1$ in the second phase.

We want to use the tutorial example to demonstrate the formulation of an MOCP with a nontrivial transition using the new \acados{} interface.
Firstly, we need to define the models for all phases, the double and the single integrator model as well as the transition model.
The single integrator model is defined as an explicit ODE using \texttt{CasADi} as follows
\begin{mypython}
	import casadi as ca
	def get_single_integrator_model() -> AcadosModel:
		model = AcadosModel()
		model.name = 'single_integrator'
		model.x = ca.SX.sym('p')
		model.u = ca.SX.sym('v')
		model.f_expl_expr = model.u
		return model
\end{mypython}
The transition model is defined as
\begin{mypython}
	def get_transition_model() -> AcadosModel:
		model = AcadosModel()
		model.name = 'transition_model'
		p = ca.SX.sym('p')
		v = ca.SX.sym('v')
		model.x = ca.vertcat(p, v)
		model.disc_dyn_expr = p
		return model
\end{mypython}
Let's assume that we already formulated the OCPs for the individual phases, the one with the single and double integrator model, which use the established Python interface of \acados, more precisely the \texttt{AcadosOcp} class.
\begin{mypython}
	def formulate_double_integrator_ocp() -> AcadosOcp:
		...
	def formulate_single_integrator_ocp() -> AcadosOcp:
		...
\end{mypython}
A multi-phase OCP can be formulated in \acados~using the \texttt{AcadosMultiPhaseOcp} class, which is created by specifying the number of stages per phase.
\begin{mypython}
	N_list = [10, 1, 15]  # 10 stages for double integrator,
	       # 1 stage for transition, 15 stages for single integrator
	ocp = AcadosMultiphaseOcp(N_list=N_list)
\end{mypython}
The novel interface allows to define the dynamics, cost and constraints of one phase utilizing the single-phase OCP formulations, i.e. the \texttt{AcadosOcp} class.
\begin{mypython}
	# define phases for single and double integrator
	phase_0 = formulate_double_integrator_ocp()
	ocp.set_phase(phase_0, 0)
	phase_2 = formulate_single_integrator_ocp()
	ocp.set_phase(phase_2, 2)
	# define the transition phase and cost
	phase_1 = AcadosOcp()
	phase_1.model = get_transition_model()
	phase_1.cost.cost_type = 'NONLINEAR_LS'
	phase_1.model.cost_y_expr = phase_1.model.x
	phase_1.cost.W = np.diag([L2_COST_P, 1e-1 * L2_COST_V])
	phase_1.cost.yref = np.array([0., 0.])
	ocp.set_phase(phase_1, 1)
\end{mypython}
Most solver options can be set in the same way as for the single-phase OCP:
\begin{mypython}
	ocp.solver_options.nlp_solver_type = 'SQP'
	ocp.solver_options.time_steps = np.array(N_list[0] * [0.4/N_list[0]]
            + [1.0] # transition stage
            + N_list[2] * [0.6 / N_list[2]])
\end{mypython}
Some additional solver options that can not vary for single-phase OCP problems can be set additionally.
\begin{mypython}
	ocp.mocp_opts.integrator_type = ['ERK', 'DISCRETE', 'ERK']
\end{mypython}
Finally, an \texttt{AcadosOcpSolver} can be created from the \texttt{AcadosMultiphaseOcp}, just as from a \texttt{AcadosOcp} object.
\begin{mypython}
	acados_ocp_solver = AcadosOcpSolver(ocp)
\end{mypython}
The interactions with the solver are independent of whether it was created from a single or multi-phase OCP formulation.

\section{Numerical experiments}
\label{sec:experiments}
This section presents three numerical case studies.
First, Section~\ref{sec:experiments_pendulum_mocp} compares MPC controllers with different control parameterizations and closed-loop costing variants on an inverted pendulum test problem.
Second, Section~\ref{sec:experiments_diff_drive_mocp} regards the task of controlling a differential drive robot directly through the voltages of actuators and compares the performance of controllers based on single- and multi-phase problem formulations.
Third, Section~\ref{sec:experiments_ptight} shows the efficiency of partial tightening with real-time iterations qualitatively replicating the benchmark results from the first partial tightening paper\cite{Zanelli2017a}.
All experiments have been carried out using \href{https://github.com/acados/acados/releases/tag/v0.3.2}{\acados{} v0.3.2} via its \texttt{Python} interface on a Laptop with an Intel i5-8365U CPU, 16 GB of RAM running Ubuntu 22.04.

Note that these experiments only compare \acados~controllers based on single- and multi-phase OCP formulations.
Comparisons with solvers based on other software packages are out of the scope of this paper.
However, since the computation times of single- and multi-phase OCP formulations within \acados~are consistent, results from existing software comparisons can be transferred to multi-phase problems if all competing solvers support MOCPs.

\subsection{Inverted pendulum on cart test problem}
\label{sec:experiments_pendulum_mocp}
We investigate the inverted pendulum on cart problem in the setting from the benchmark presented by Frey et al.~\cite{Frey2024b} and compare different controllers with different closed-loop costing variants and control parametrizations.
The code of the orignal benchmark has been adapted to incorporate the new controllers\footnote{\url{https://github.com/FreyJo/GNRK_benchmark}}.

The differential state of the model is $x = [p, \theta, s, \omega]^\top$ with cart position $p$, cart velocity $s$, angle of the pendulum $\theta$ and angular velocity $\omega$.
The control input $v$ is a force acting on the cart in the horizontal plane.
The system dynamics can be found e.g. in~\cite{Verschueren2021}.
In our OCP formulation, $v$ is constrained to be in $[-40, 40]$.
The simulation starts with an initial state $\bar x_0 = [0, \frac{\pi}{5}, 0, 0]^\top$.
The goal is to drive all states to zero, i.e. the unstable upright position.
We formulate the following nonlinear least-squares cost consisting of quadratic costs on states and controls and a term penalizing a position $p$ outside of $[-1, 1] =: [p_{\min}, p_{\max}]$, namely
\begin{align*}
    l(x, u) &= x^\top Q x + v^\top R v + \gamma \cdot (\max(p_{\min} \shortminus p, 0))^2 + \gamma \cdot (\max(p \shortminus p_{\max}, 0)^2 \nonumber
\end{align*}
where
the cost weights are chosen as $\gamma=5\cdot 10^4$, $Q = \mathrm{diag}(100, 10^3, 0.01, 0.01) $, $ R = 0.2 $.
The terminal cost term is set to $M\ind{pend}(x) = x^\top P x$, where $P$ is obtained as solution of the discrete algebraic Riccati equation with cost and dynamics linearized at the steady state.

\begin{figure}
    \vspace{-.0cm}
    \centering
    	\newcommand{\startTime}{0}
	\newcommand{\numIntervalsFirstRange}{4}
	\newcommand{\numIntervalsSecondRange}{5}
	\newcommand{\subintervalLengthZero}{0.02}
	\newcommand{\timeLengthFirstRange}{0.28}
	\newcommand{\timeLengthSecondRange}{3.7}

	\pgfmathsetmacro{\subintervalLengthFirstRange}{(\timeLengthFirstRange) / \numIntervalsFirstRange}
	\pgfmathsetmacro{\subintervalLengthSecondRange}{\timeLengthSecondRange / \numIntervalsSecondRange}
	\pgfmathsetmacro{\Ttotal}{\timeLengthSecondRange + \timeLengthFirstRange + \subintervalLengthZero}

	\pgfmathsetmacro{\Tone}{\timeLengthFirstRange + \subintervalLengthZero}
	\pgfmathsetmacro{\numIntervals}{1 + \numIntervalsFirstRange + \numIntervalsSecondRange}

	\pgfmathsetmacro{\tikheight}{0.05}
	\begin{tikzpicture}[scale=1.5]
		\node[above] at (\startTime - 0.7,-\tikheight) {Grid A};

		\pgfmathsetmacro{\totalTimeLength}{\timeLengthFirstRange + \timeLengthSecondRange}

		\draw[->] (\startTime,0) -- ({\startTime + 1.1 * \totalTimeLength},0) node[right] {$t$};

		\draw (\startTime, \tikheight) -- (\startTime, -\tikheight);
		\draw ({\startTime + \subintervalLengthZero}, \tikheight) -- ({\startTime + \subintervalLengthZero}, -\tikheight);

		\foreach \i in {1,...,\numIntervalsFirstRange} {
		  \draw ({\startTime + \subintervalLengthZero + (\i-1)*\subintervalLengthFirstRange}, \tikheight) -- ({\startTime + \subintervalLengthZero + (\i-1)*\subintervalLengthFirstRange}, -\tikheight);
		}

		\foreach \i in {0,...,\numIntervalsSecondRange} {
		  \draw ({\startTime + \subintervalLengthZero + \timeLengthFirstRange + \i*\subintervalLengthSecondRange}, \tikheight) -- ({\startTime + \subintervalLengthZero + \timeLengthFirstRange + \i*\subintervalLengthSecondRange}, -\tikheight);
		}

		\node[below] at (\startTime,-\tikheight) {$\startTime$};
		\node[below] at ({\startTime + \totalTimeLength},-\tikheight) {$T=\Ttotal$};

		\node[below] at (\startTime+\subintervalLengthZero-.1,+0.4) {$\dtplant$};
		\draw[->] (\startTime+\subintervalLengthZero / 2, 0.18) -> (\startTime+\subintervalLengthZero / 2, 0.1);

		\node[below] at (\Tone+0.05, +0.55) {\footnotesize{$ T_1 = \Tone $}};
		\draw[->] (\Tone, 0.34) -> (\Tone, 0.15);

	  \end{tikzpicture}
      \pgfmathsetmacro{\dtlong}{(\timeLengthFirstRange + \timeLengthSecondRange) / \numIntervals}

      \begin{tikzpicture}[scale=1.5]
		\node[above] at (\startTime - 0.7,-\tikheight) {Grid B};

		\pgfmathsetmacro{\totalTimeLength}{\timeLengthFirstRange + \timeLengthSecondRange}

		\draw[->] (\startTime,0) -- ({\startTime + 1.1 * \totalTimeLength},0) node[right] {$t$};

		\draw (\startTime, \tikheight) -- (\startTime, -\tikheight);
		\draw ({\startTime + \subintervalLengthZero}, \tikheight) -- ({\startTime + \subintervalLengthZero}, -\tikheight);

		\foreach \i in {0,...,\numIntervals} {
		  \draw ({\startTime + \subintervalLengthZero + (\i)* \dtlong}, \tikheight) -- ({\startTime + \subintervalLengthZero + (\i)*\dtlong}, -\tikheight);
		}

		\node[below] at (\startTime,-\tikheight) {$\startTime$};
		\node[below] at ({\startTime + \totalTimeLength},-\tikheight) {$T=\Ttotal$};

		\node[below] at (\startTime+\subintervalLengthZero-.1,+0.4) {$\dtplant$};
		\draw[->] (\startTime+\subintervalLengthZero / 2, 0.18) -> (\startTime+\subintervalLengthZero / 2, 0.1);

	  \end{tikzpicture}
    \caption{Two time grids considered in the benchmark. Both grids start with an interval of length $\Delta t = 0.02 $.
    The grids are visualized for $N = 10$ shooting intervals.
    For Grid A, the interval $[\Delta t, T_1] $ is divided into $N_1 = 4$ equidistant intervals, and $[T_1, T]$ into $N_2 = N-N_1-1 =5$ equidistant intervals.
    For Grid B, the interval $[\Delta t, T] $ is divided into $N-1$ intervals.
    \label{fig:timegrids}
    }
\end{figure}
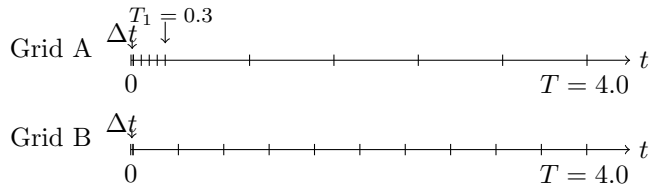

The controllers use two different discretization grids, visualized in Figure~\ref{fig:timegrids}.
The plant is simulated with a time step of $\dtplant = 0.02 \mathrm{s}$.
All controllers use the GNRK cost discretization described in our recent paper\cite{Frey2024b} and Section \ref{sec:cost_integration} with a Gauss-Radau IIA method of order 7 on each shooting interval.
Additionally, we compare various control parameterizations, namely piecewise polynomials of different degree on all but the first shooting interval $[0, \Delta t]$, see Section~\ref{sec:pw_polynomial}, and closed-loop costing formulation, see Section~\ref{sec:clc}.
In order to conveniently formulate the models with piecewise polynomial controls, we added the functionality \texttt{AcadosModel.reformulate\_with\_polynomial\_control()}.

\begin{figure}
    \centering
    \includegraphics[width=.7\columnwidth]{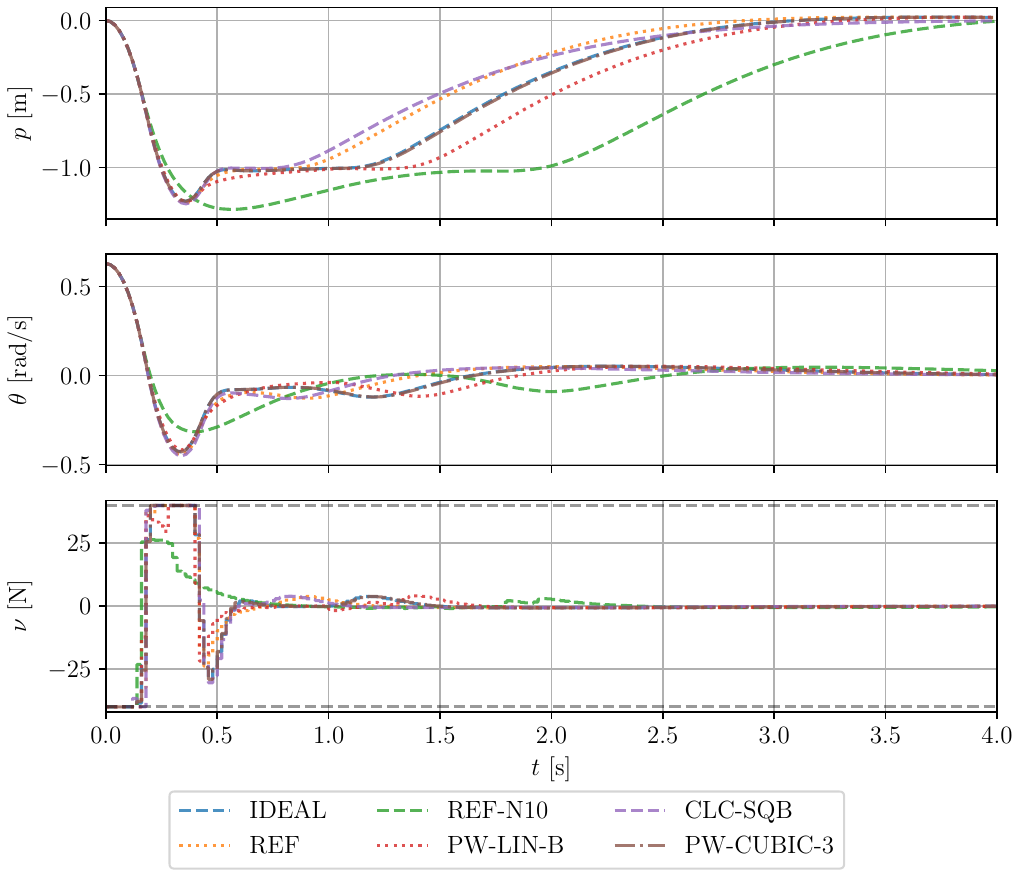}
    \caption{Closed-loop trajectories corresponding to the benchmark results in Table~\ref{tab:pendulum_closed_loop_costing}.
    Values for the postion $p$ outside of $[-1, 1]$ are heavily penalized and result in high suboptimality, which is reported in Figure~\ref{fig:pendulum_pareto_multiphase} and Table~\ref{tab:pendulum_closed_loop_costing}.
    \label{fig:pendulum_clc_traj}
    }
\end{figure}

\begin{figure}
	\centering
	\includegraphics[width=.7\columnwidth]{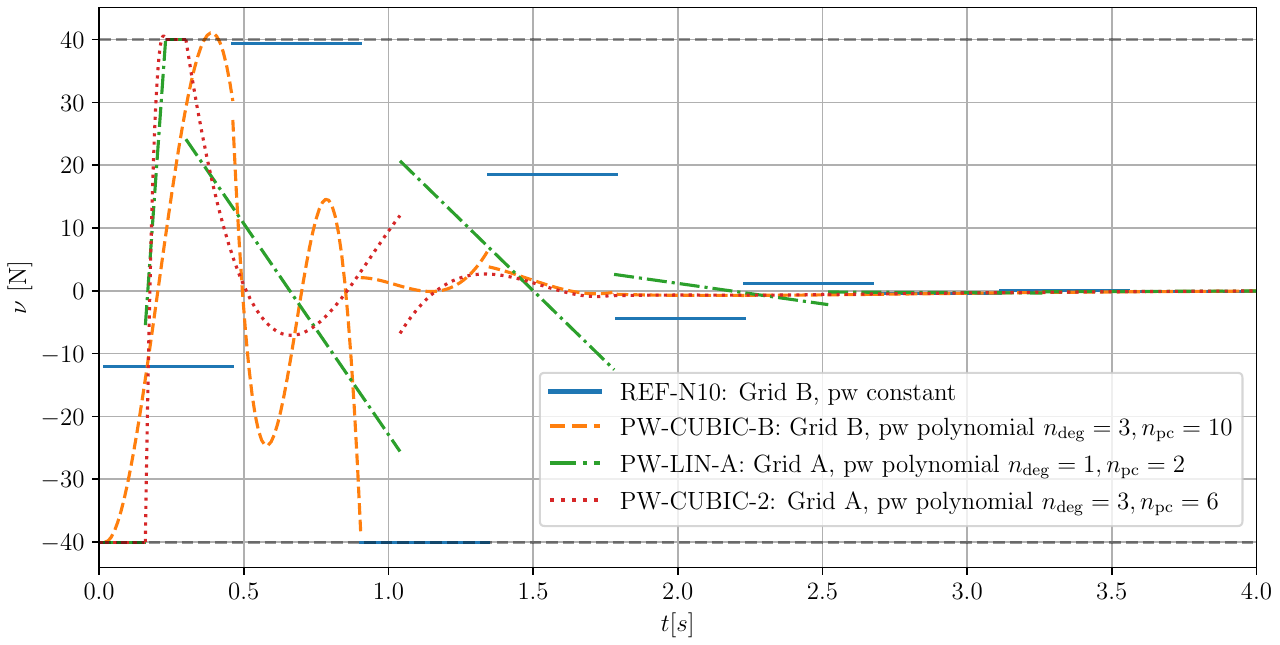}
	\caption{Open-loop control trajectories corresponding to the first problem solved in the closed-loop simulation visualized in Figure~\ref{fig:pendulum_clc_traj}. Details on the controller variants are in Table~\ref{tab:pendulum_closed_loop_costing}.
		\label{fig:pendulum-polynomial-u}
	}
\end{figure}

An overview of all controller variants considered is presented in Table~\ref{tab:pendulum_closed_loop_costing}, some of the closed-loop trajectories are plotted in Figure~\ref{fig:pendulum_clc_traj} and some of the open-loop control trajectories corresponding to the first problem are shown in Figure~\ref{fig:pendulum-polynomial-u}.

\begin{figure}
	\centering
	\includegraphics[width=.7\columnwidth]{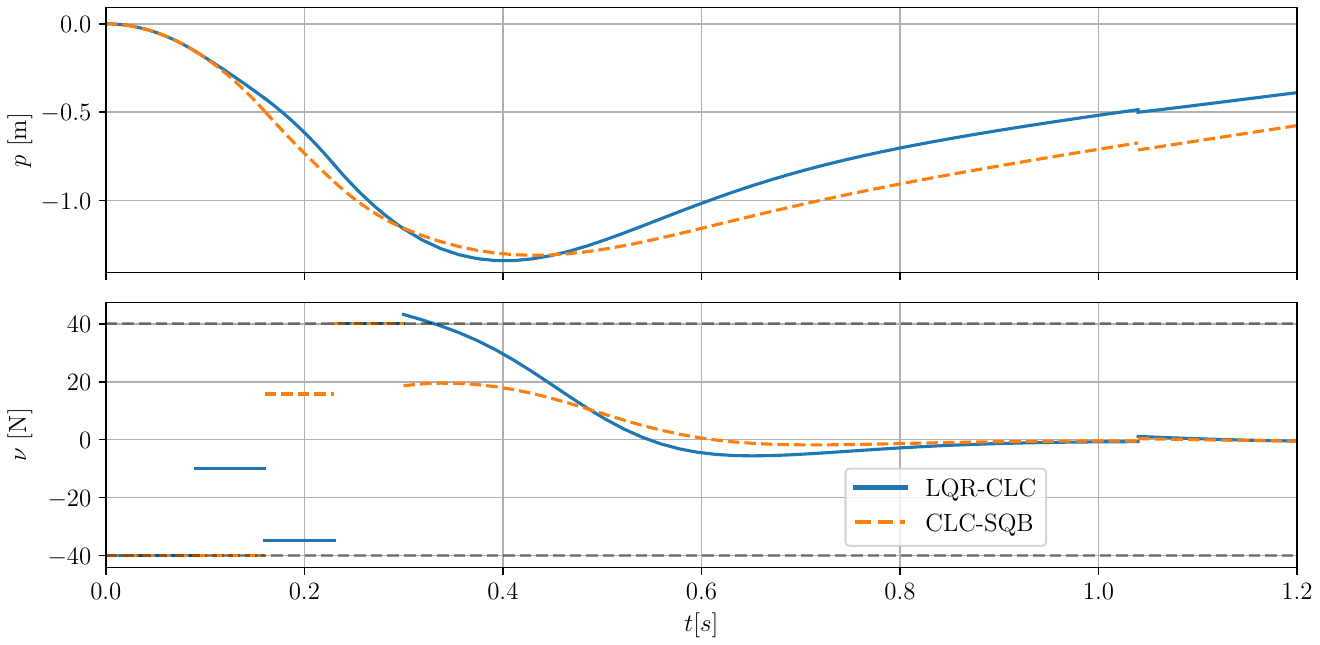}
	\caption{Multiple-shooting with closed-loop costing. The solver iterates after performing 3 iterations on the first problem in the closed-loop scenario in Figure~\ref{fig:pendulum_clc_traj} are visualized.
		\label{fig:pendulum-clc-ms}
	}
\end{figure}

\begin{table*}
\centering
\caption{Closed-loop comparison of different controller variants on the pendulum test problem with different discretization grids and the control parametrizations, including closed-loop-costing and piecewise polynomial controls.
}
\label{tab:pendulum_closed_loop_costing}
{\footnotesize
\begin{tabular}{lrrcrrrr}
\toprule
variant ID & $N$ & Grid & control parametrization & comp. time / iter [ms] & rel. suboptimality [\%] \\ \midrule
IDEAL &200 & B &pw. constant &3.69 &0.00\\
REF &20 & B &pw. constant &0.43 &3.72\\
REF-N10 &10 & B &pw. constant &0.21 &267.57\\
PW-LIN-B &10 & B &pw. polynomials, $n\ind{deg} =  1, n\ind{pc} = 2 $ &0.23 &13.54\\
PW-CUBIC-B &10 & B &pw. polynomials, $n\ind{deg} =  3, n\ind{pc} = 10 $ &0.24 &13.59\\
\midrule
LQR-CLC &10 & A & unconstrained LQR CLC &0.22 &3.41\\
CLC-SQB &10 & A & Squashed  + prog. barrier CLC &0.23 &5.13\\
PW-CONST-A &10 & A &pw. constant &0.22 &0.49\\
\midrule
PW-LIN-A &10 & A &pw. polynomials, $n\ind{deg} =  1, n\ind{pc} = 2 $ &0.23 &0.40\\
PW-QUAD-1 &10 & A &pw. polynomials, $n\ind{deg} =  2, n\ind{pc} = 4 $ &0.24 &0.06\\
PW-QUAD-2 &10 & A &pw. polynomials, $n\ind{deg} =  2, n\ind{pc} = 10 $ &0.26 &0.15\\
\midrule
PW-CUBIC-1 &10 & A &pw. polynomials, $n\ind{deg} =  3, n\ind{pc} = 4 $ &0.23 &0.07\\
PW-CUBIC-2 &10 & A &pw. polynomials, $n\ind{deg} =  3, n\ind{pc} = 6 $ &0.24 &0.03\\
PW-CUBIC-3 &10 & A &pw. polynomials, $n\ind{deg} =  3, n\ind{pc} = 10 $ &0.26 &0.02\\
\bottomrule
\end{tabular}
}
\end{table*}

\begin{figure}
	\centering
	\includegraphics[width=.7\columnwidth]{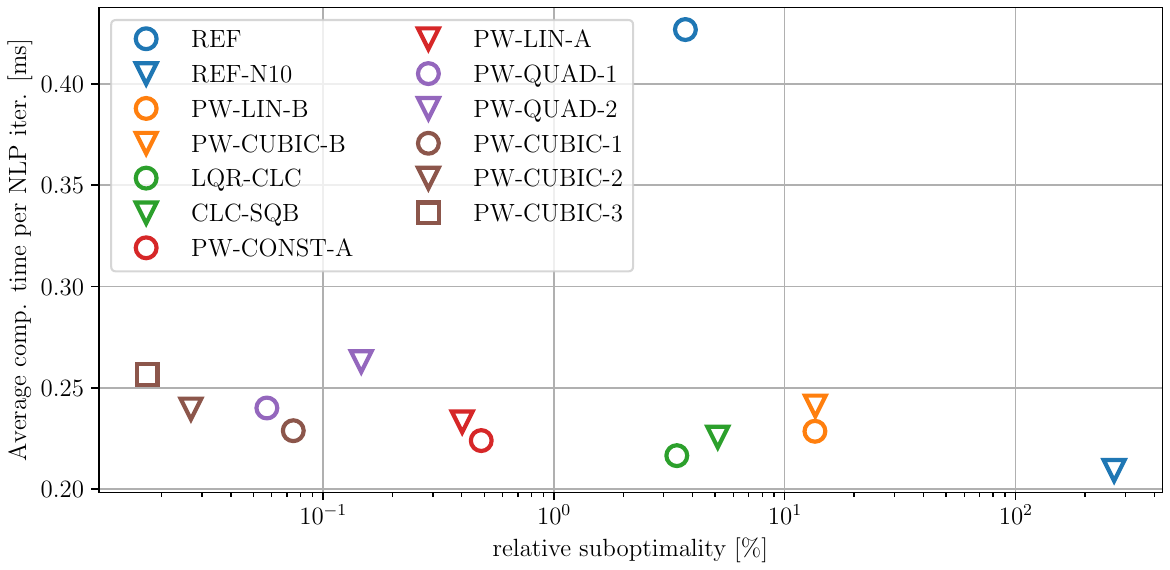}
	\caption{Pareto plot comparing the controllers in Table~\ref{tab:pendulum_closed_loop_costing} in terms of average computation time per NLP solver iteration and relative suboptimality with respect to the controller IDEAL which is not included in this Figure.
		\label{fig:pendulum_pareto_multiphase}
	}
\end{figure}

\subsubsection{Controller variants}
The controller labeled IDEAL uses a uniform time discretization corresponding to time steps $\dtplant$ with piecewise constant controls, i.e. no mismatch between the plant and the OCP model.
This controller gives the best closed-loop performance, but has the highest computational complexity.
The controller REF corresponds to the variant which gave the best tradeoff between computation time and closed-loop performance in the original benchmark\cite{Frey2024b}.
The goal is to derive a discrete-time OCP variant which uses less shooting intervals than this controller in order to reduce the associated computational complexity, while attaining similar closed-loop performance.

The variant REF-N10 uses the same kind of grid, but divides the interval $[\Delta t, T]$ into $9$ instead of $19$ shooting intervals.
Thus, the computational complexity associated with one iteration is roughly halved.
However, the closed-loop performance is highly degraded.
The controller variant PW-LIN-B has only one difference, which is to use a piecewise linear control parametrization on all but the first shooting interval.
This results in a much better performance compared to REF-N10, with marginally higher computational burden, indicating that the cost-to-go approximation on this coarse grid is limited by the control parameterization.

We evaluate two closed-loop costing (CLC) variants.
Both CLC controllers use Grid A, where the second part of the time horizon, $[T_1, T]$ is implemented with a second phase, see Section~\ref{sec:clc}.
For the LQR-CLC controller, the control law is the stabilizing LQR control law $ \kappa\ind{LQR}(x) = Kx$.
The control input constraints are not enforced on the CLC horizon.
This controller results in a closed-loop performance that is similar to the one of controller REF proposed in the original benchmark\cite{Frey2024b} at half the compuation time.

The other CLC controller is labeled CLC-SQB and uses the same linear feedback law with additional squashing to impose the given control limits, $ \kappa\ind{squash}(x) = \sigma(Kx) $.
Additionally, a barrier term $\beta(\sigma(Kx))$ is added on the CLC horizon.
The resulting controller can therefore not plan with violations of the control bounds in the CLC horizon.
This controller yields slightly worse closed-loop performance compared to the standard LQR-CLC controller.
We attribute this observation to the fact that the CLC controller with squashing naturally has a more conservative plan, which results in some additional suboptimality.
Figure~\ref{fig:pendulum-clc-ms} visualizes intermediate iterates of the controllers LQR-CLC and CLC-SQB and gives some insights into CLC.
It can be seen that the controls on the CLC horizon are continuous, since $\kappa$ and the state trajectory are continuous.
However, for the intermediate iterates of the multiple-shooting formulation, the shooting gaps, which are visible in Fig. \ref{fig:pendulum-clc-ms} result in gaps in the controls within the CLC horizon which are closed at convergence.
Lastly, we observe that the LQR-CLC controller violates the control bounds in its plan, while the squashed variant CLC-SQB inherently respects them.

As an alternative to the CLC controllers, some variants are included, which use discretization Grid A, but a standard control horizon on the second part of the grid $[T_1, T]$ and piecewise polynomial controls on $[\Delta t, T]$ with the control bounds enforced on $n\ind{pc}$ equidistant points on every shooting interval, see Section~\ref{sec:pw_polynomial}.

The controllers PW-CONST-A and PW-LIN-A exactly enforce the control bounds, are among the controllers with the lowest computational complexities and result in a relative suboptimality below 0.5 $\%$.
In particular, they are able to roughly halve the computational complexity of REF, the best controller in the original benchmark\cite{Frey2024b}, while reducing the suboptimality by a factor greater than 7.

\subsubsection{Discussion}
The fact that controller PW-CONST-A has a lower computational burden compared to the CLC variant is attributed to the fact that the control law $\kappa(\cdot)$ and its derivatives have to be evaluated often within every step and Newton iteration of the implicit integrator.
The additional computational burden of a larger input dimension, which is 0 or 1 on the latter part of the horizon, is rather small for efficient state-of-the-art QP solvers.

There are certainly other examples in which a CLC controller is computationally more attractive.
In addition to the test setting, this of course depends on the underlying software framework for linearization and solution of the subproblems.
In particular, for modern OCP-structure exploiting software frameworks, the computational burden associated with an extra control variable is less significant compared to when CLC was introduced in the early 1990s \cite{Kouzoupis2018}.

A similar comparison of CLC controllers with respect to obvious alternatives is not to be found in the existing literature on CLC\cite{Nicolao1998, DeNicolao1996a, Diehl2004f, Magni2004, Quirynen2014a}.
In addition to the LQR-CLC controller described above, we compare a related variant, which uses $N_2=1$, i.e. single shooting on the CLC horizon, and instead performs five steps of IRK on this interval.
This controller variant does not converge in the first simulation step of the scenario and has a computational complexity similar to LQR-CLC.
This indicates that the multiple shooting CLC implementation proposed in this paper has desirable properties compared to the single-shooting CLC variants used in previous works.

Overall, the comparison of CLC with respect to nonuniform OCP discretizations shows that the latter are very competitive.
On this test example, the nonuniform grid variant outperforms the LQR closed-loop costing based controller in terms of closed-loop performance by a factor of 7 while requiring an equal amount of computational resources.

Finally, we regard controllers with piecewise polynomial control parameterizations of degree $n\ind{deg}>1$ which enforce the control bounds on $n\ind{pc}$ intermediate points.
The open-loop solutions corresponding to the first problems within the closed-loop simulation are visualized for a few controllers in Figure~\ref{fig:pendulum-polynomial-u}.
We observe small violations of the open-loop control trajectories with respect to the control bounds, e.g. at $t \approx 0.2$ for PW-CUBIC-2 and $t\approx 0.35$ for PW-CUBIC-B.

All variants in Figure~\ref{fig:pendulum-polynomial-u} use 10 shooting intervals.
Comparing REF-N10 in Figure~\ref{fig:pendulum-polynomial-u} with the optimal closed-loop control trajectory in Figure~\ref{fig:pendulum_clc_traj}, we observe that this planned controls are qualitatively very different from the optimal ones, resulting in a bad cost-to-go approximation and thus closed-loop performance.
In contrast, PW-CUBIC-2 can capture the qualitative behavior of the optimal closed-loop control trajectory well.
In particular, the relatively short intervals in the first part of the horizon enable the cubic polynomials to approximate the optimal bang-bang solution in this phase well.
In contrast, the variant PW-CUBIC-B, which uses the same control parametrization as PW-CUBIC-2, but on grid A, approximates this optimal bang-bang behavior with a slower transition, see Figure~\ref{fig:pendulum_clc_traj}, resulting in significantly higher closed-loop cost, see Table~\ref{tab:pendulum_closed_loop_costing}.
A visual comparison of the open-loop trajectories of PW-CUBIC-2 and PW-LIN-A shows that the optimal bang-bang behavior in the first part of the horizon is approximated similarly, while on the first long interval, starting at $t=0.3$, the cubic polynomials capture the optimal behavior qualitatively better, as piecewise linear functions do not have an inflection point.

When comparing controllers with the same parametrization, but different grids on which the control bounds are enforced, i.e. varying $n\ind{pc}$, we can observe two effects.
Firstly, comparing PW-QUAD-1 and PW-QUAD-2, we see that enforcing control bounds on a finer grid can result in worse closed loop performance.
This can be attributed to the fact that enforcing the control bounds everywhere using polynomials results in a more conservative approximation of the cost-to-go, since the planned controls are approximated within the bounds.
On the other hand, regarding PW-CUBIC-1, PW-CUBIC-2, PW-CUBIC-3, we observe that enforcing constraints on a tighter grid can improve the approximation quality.

Overall, the general trend in our experiments shows that additional degrees of freedom in piecewise polynomial control parametrizations can result in a strongly improved closed-loop performance, if the control bounds are enforced on a sufficiently fine grid, which we attribute to a better approximation of the closed-loop cost.
In particular, PW-CUBIC-3 results in better performance than PW-QUAD-2, which in turn improves on PW-LIN-A.

Of course one needs to be careful, when it comes to drawing general conclusions from the presented results on piecewise polynomial control parameterizations.
The additional computational cost of polynomial control parameterizations of higher degrees and handling the corresponding bounds, will be much more significant, with growing $n_\contcontrl$.
In addition to the dimensions, the optimal choice of control parameterization and where to enforce control bounds depends on the qualitative behavior of the continuous-time optimal solution (bang-bang or continuous), the time discretization grid, the available solvers and computational resources.
However, the possibility of formulating such problems within an efficient software package can significantly improve the closed-loop performance of MPC controllers.

\subsection{Differential drive robot with actuator model and economic cost}
\label{sec:experiments_diff_drive_mocp}
We consider a differential drive robot and develop an NMPC controller which takes the underlying actuators into account allowing it to consider their power consumption in the cost function.
The code to reproduce the results presented in this section is publicly available\footnote{\url{https://github.com/FreyJo/ocp_solver_benchmark/blob/main/experiments/actuator_diff_drive.py}}

\subsubsection{Modelling}
This subsection describes the differential drive model with and without actuators~\cite{Dhaouadi2013}.
For simplicity, the dynamic model which disregards the actuators is presented first.
This model consists of the state vector $x\ind{simple} = [p\ind{x}, p\ind{y}, v, \theta, \omega]$, where $ p\ind{x}, p\ind{y} $ denotes the robots position in x- and y-coordinates, $v$ the velocity of the robot, $\theta$ the heading angle, and $\omega$ the angular velocity.
The input of this model are $u\ind{simple} = [\tau\ind{r}, \tau\ind{l}] $, the torques applied to the right and left driving wheel, respectively, in $\mathrm{Nm}$.
The evolution of the system is described by the ODE
\begin{align}
\dot{x}\ind{simple} =
\begin{bmatrix}
v \cos \theta \\
v \sin \theta \\
\frac{a_1 + m\ind{c} d \omega^2}{m + a_2} \\
\frac{L a_3 - m\ind{c} d \omega v}{I + L^2 a_2}
\end{bmatrix},
\label{eq:simple_diff_drive}
\end{align}
where the shorthands $a_1 = \frac{\tau\ind{r} + \tau\ind{l}}{2}$, $a_2 = \frac{2 I\ind{w}}{R^2}$, $a_3 = \frac{\tau_r - \tau_l}{R}$ are introduced.
The model contains the following parameters, which we assume to be constant:
The total mass of the robot $m = 220 \mathrm{kg}$, the mass of the robot without the wheels and the rotating parts of the actuators $m\ind{c} = 200 \mathrm{kg}$, the robot's moment of inertia around the center of mass $I=9.6 \, \mathrm{kg} \cdot \mathrm{m}^2$, the combined moment of inertia about the wheel’s axis of a driving wheel and the rotating part of the actuator $I\ind{w} = 0.1 \mathrm{kg} \cdot \mathrm{m}^2$.

In addition, we consider a more accurate model of the robot which takes the two actuators, namely a motor on each driving wheel into account.
The actuator model consists of the state $x\ind{act} = [p\ind{x}, p\ind{y}, v, \theta, \omega, I\ind{r}, I\ind{l}]$, where $ I\ind{r}, I\ind{l}$ are the currents in the motor driving the right and left wheel, respectively, and the control input $u\ind{act} = [V\ind{r}, V\ind{l}]$, where $V\ind{r}, V\ind{l}$ denote the voltages applied to the motors.
The evolution of the states common for both models is given by \eqref{eq:simple_diff_drive}, where the motor torques $\tau\ind{r}, \tau\ind{l}$ are substituted by $\tau\ind{r} = K_1 I\ind{r}$ and $\tau\ind{l} = K_2 I\ind{l}$.
In addition, the dynamics of the accurate model are given by
\begin{align}
\begin{bmatrix}
\dot{I}\ind{r} \\
\dot{I}\ind{l}
\end{bmatrix}
=
\begin{bmatrix}
- \frac{K_1 \psi_1 - R\ind{act} I\ind{r} + V\ind{r}}{L\ind{act}} \\
- \frac{K_2 \psi_2 - R\ind{act} I\ind{l} + V\ind{l}}{L\ind{act}}
\end{bmatrix}
\end{align}
with
\begin{align}
\psi_1 & = \frac{\dot{p}\ind{x} \cos \theta + \dot{p}\ind{y} \sin \theta + L \omega}{R}, \\
\psi_2 & = \frac{\dot{p}\ind{x} \cos \theta + \dot{p}\ind{y} \sin \theta - L \omega}{R}.
\end{align}
This accurate model additionally contains the motor constants $ K_1 = 1.0 $, $K_2 = 1.0$, the coil inductance $L\ind{act} = 10^{-4} \mathrm{H} $, the coil resistance $R\ind{act} = 0.05 \Omega$ as parameters which we assume to be constant.

\subsubsection{Optimal control problem formulation}
We want to minimize the following cost term
\begin{align}
l\ind{act}(x\ind{act}, u\ind{act}) = 
x\ind{act}\transp Q x\ind{act} + \abs{V\ind{r} I\ind{r}} + \abs{V\ind{l} I\ind{l}}
\label{eq:ddrive_cost}
\end{align}
with $Q = \mathrm{diag}(10^3, 10^3, 10^{-4}, 1, 10^{-3}, 0.5, 0.5) $.
Note that the second summand $\abs{V\ind{r} I\ind{r}} + \abs{V\ind{l} I\ind{l}}$ is an economic cost term corresponding to the power consumption.
It is implemented in~\acados~by introducing slack variables $s\ind{low}$, $s\ind{up}$ and constraints $s\ind{low} \leq [V\ind{r} I\ind{r}, V\ind{l} I\ind{l}]\transp \leq s\ind{up}$ and replacing the term $\abs{V\ind{r} I\ind{r}} + \abs{V\ind{l} I\ind{l}}$ in the cost with $s\ind{low, 1} + s\ind{low, 2} + s\ind{up, 1} + s\ind{up, 2}$.
When using the actuator model, for the full horizon, we use the terminal cost term:
\begin{align}
E\ind{act}(x\ind{act}) = x\ind{act}\transp Q x\ind{act}
\end{align}

We want to investigate approximate MOCP formulations which use the simple differential drive model on the latter part of the horizon.
The cost term in \eqref{eq:ddrive_cost} can be approximated using the simple model by disregarding the economic cost term, i.e. using
\begin{align}
	l\ind{simple}(x\ind{simple}, u\ind{simple}) =
	x\ind{simple}\transp \tilde{Q} x\ind{simple} + u\ind{simple}\transp \tilde{R} u\ind{simple}
	\label{eq:ddrive_cost_approx}
\end{align}
with $\tilde{Q} = \mathrm{diag}(10^3, 10^3, 10^{-4}, 1, 10^{-3}) $, $\tilde{R} = \mathrm{diag}(0.5, 0.5) $.
As a terminal cost of the first phase with the actuator model, cf. \eqref{eq:dmocp}, we define
\begin{align}
	E\ind{trans}(x\ind{act}) = 0.5 \Delta t\ind{trans} \cdot (I\ind{r}^2 + I\ind{l}^2),
\end{align}
where $\Delta t\ind{trans}$ denotes the length of the last shooting interval on which the actuator model is used.
This cost term corresponds to the tracking part of the cost in \eqref{eq:ddrive_cost} and is needed since $I\ind{r}$, $I\ind{l}$ are not included in the simple model.
When using the simple model at the end of the horizon, we use the terminal cost term:
\begin{align}
	E\ind{simple}(x\ind{simple}) = x\ind{simple}\transp \tilde{Q} x\ind{simple}.
\end{align}
For both models, we impose the following path constraints on the state
\begin{align}
0 &\leq v \leq 1, \\
-0.5 &\leq \omega \leq 0.5.
\end{align}
Additionally, we impose the following constraints on the control inputs when using the actuator model
\begin{align}
 -10 \leq V\ind{l} \leq 10,\\
 -10 \leq V\ind{r} \leq 10.
\end{align}
Respectively, for the simple model, we impose
\begin{align}
-60 \leq \tau\ind{r} \leq 60,\\
-60 \leq \tau\ind{l} \leq 60.
\end{align}

\begin{figure}[]
	\centering
	\includegraphics[width=.7\columnwidth]{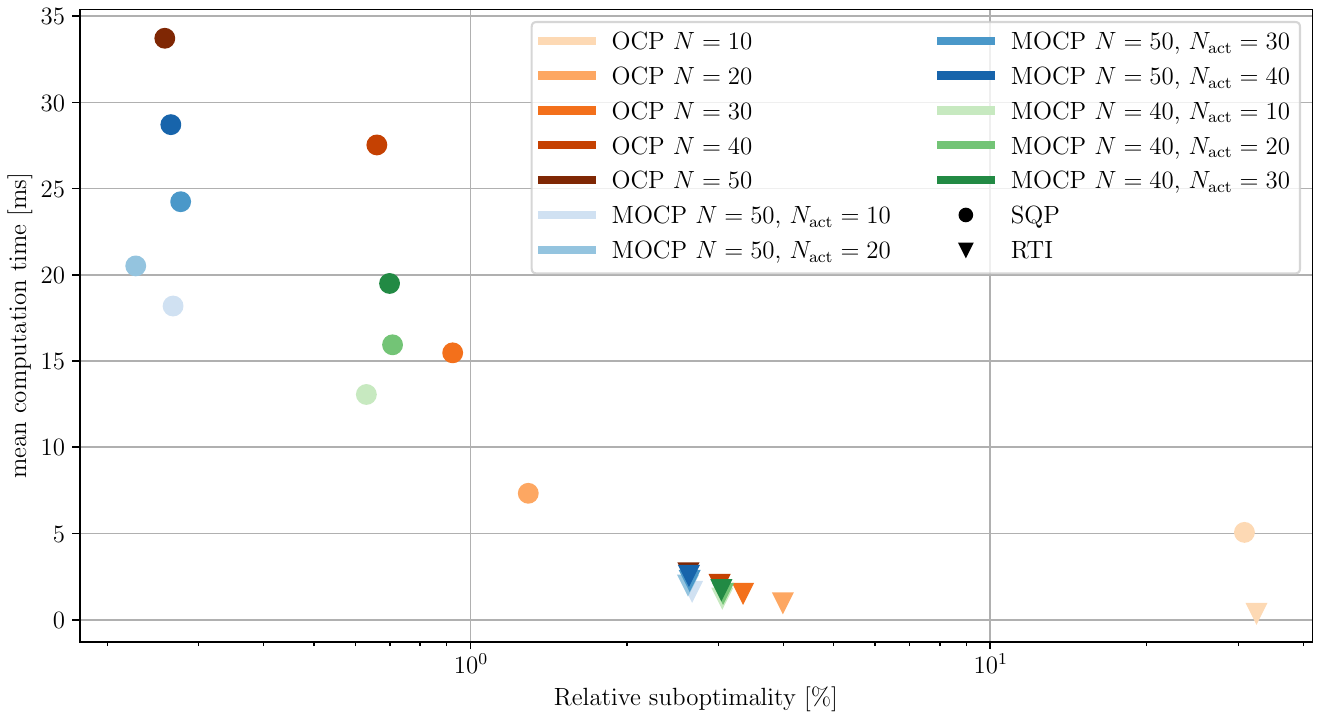}
	\caption{Pareto plot comparing different controller variants for directly controlling the actuators of a differential drive robot with an economic cost.
		All MOCP variants use $N=50$ shooting intervals and $N-N\ind{act}$ of an approximate model which disregards the actuator dynamics.
		\label{fig:diff_drive_pareto}
	}
\end{figure}
\subsubsection{Comparing controller variants}
We compare different controller variants.
All variants use discretization steps of $\Delta t = 0.1 \mathrm{s}$, the model dynamics are integrated using an implicit Runge-Kutta method with a Legendre Butcher tableau of order 6 on each shooting interval.
The controllers are based on an \acados~full-step SQP solver using a Gauss-Newton Hessian approximation and solving the QPs with \texttt{HPIPM} without condensing.

We consider MOCP variants with $N=40$ and $N=50$ shooting intervals and vary $N\ind{act}$, the number of shooting intervals on which the actuator model is used, afterwards, we use a transition stage and $N-N\ind{act}$ shooting intervals with the simple model.
For the single-phase OCP variants, we vary the number of shooting intervals $N$, as this formulation does not allow for a transition.

We compare several controllers in a closed-loop simulation of $20\mathrm{s}$ using the actuator model starting at the initial state $x_0 = (1, 1, 0, \pi, 0, 0, 0) $, exactly integrating the cost in \eqref{eq:ddrive_cost} using 10 integrator steps for a sampling time of $ 0.1 \mathrm{s}$.
In order to evaluate relative suboptimality, a reference controller with $N=60$ shooting intervals of the actuator model is used.
Controllers with fully converged SQP and RTI variants, which only solve a single QP subproblem are evaluated in the same plot.

The Pareto plot in Figure~\ref{fig:diff_drive_pareto} visualizes the suboptimality and computation time of the various variants.
We observe that the Pareto front is dominated by the MOCP variants.
Note that the variant OCP $N=50$ correspond to the controller MOCP $N = 50, N\ind{act} = 50$.
The solver variants OCP $N=30$ and MOCP $N=50 , N\ind{act}=10$ with SQP require a similar computation time, however the MOCP variant results in a three fold lower relative suboptimality.
On the other hand, MOCP with $N=50, N\ind{act}=10$ delivers a similar suboptimality compared to the variant OCP $N=50$, while only requiring roughly half of the computation time.
Overall, the MOCP variants dominate the Pareto front in large parts and it is of course possible to generate many more combinations using MOCP formulations or by varying the number of SQP iterations.
These results show that deriving an approximate OCP formulation and using it within an MOCP can result in an NMPC controller which is able to outperform NMPC controllers that only use a single-phase OCP formulation.

\FloatBarrier
\subsection{Partial tightening} \label{sec:experiments_ptight}
\begin{table}
	\centering
	\caption{Performance overview of different controllers with RTI and partial tightening. Maximum timings are over the closed-loop simulation in Figure~\ref{fig:zanelli_ptight_trajectories} and are given in [ms] and relative suboptimality is evaluated by comparing with the controller variant $N=100, N_{\mathrm{exact}}= 100$.\label{tab:ptight_performance}}
	\begin{tabular}{lcccccc}
		\toprule
		& \multicolumn{3}{c}{computation times} & \\
		Variant & preparation & feedback & total & rel. subopt. \%\\ \midrule
		$N=100, N_{\mathrm{exact}}= 5$      & 0.23 & 0.32 & 0.54 & 21.25\\
		$N=100, N_{\mathrm{exact}}= 10$     & 0.25 & 0.36 & 0.59 & 16.61\\
		$N=100, N_{\mathrm{exact}}= 20$     & 0.24 & 0.44 & 0.65 & 12.25\\
		$N=100, N_{\mathrm{exact}}= 50$     & 0.25 & 1.15 & 1.36 & 5.15\\
		$N=100, N_{\mathrm{exact}}= 100$    & 0.23 & 2.87 & 3.05 & 0.00\\
		$N=50, N_{\mathrm{exact}}= 50$      & 0.10 & 2.17 & 2.27 & 990.43\\
		\bottomrule
	\end{tabular}
\end{table}

\begin{figure}[]
	\centering
	\includegraphics[width=.7\columnwidth]{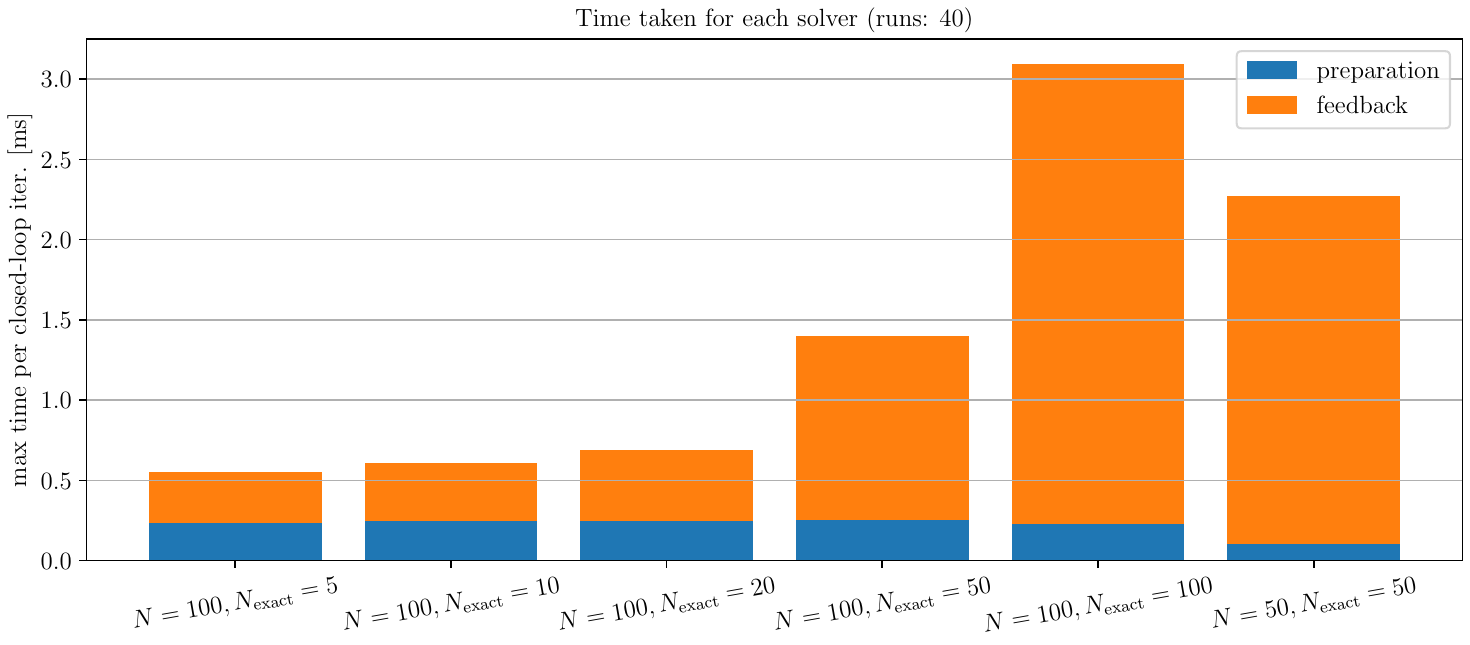}
	\caption{Maximum timings for the preparation and feedback phase over the closed-loop simulation  shown in Figure~\ref{fig:zanelli_ptight_trajectories}.
	The timings are also given in Table~\ref{tab:ptight_performance}.
	\label{fig:zanelli_ptight_rti_timings}
	}
\end{figure}

\begin{figure}[]
	\centering
	\includegraphics[width=.9\columnwidth]{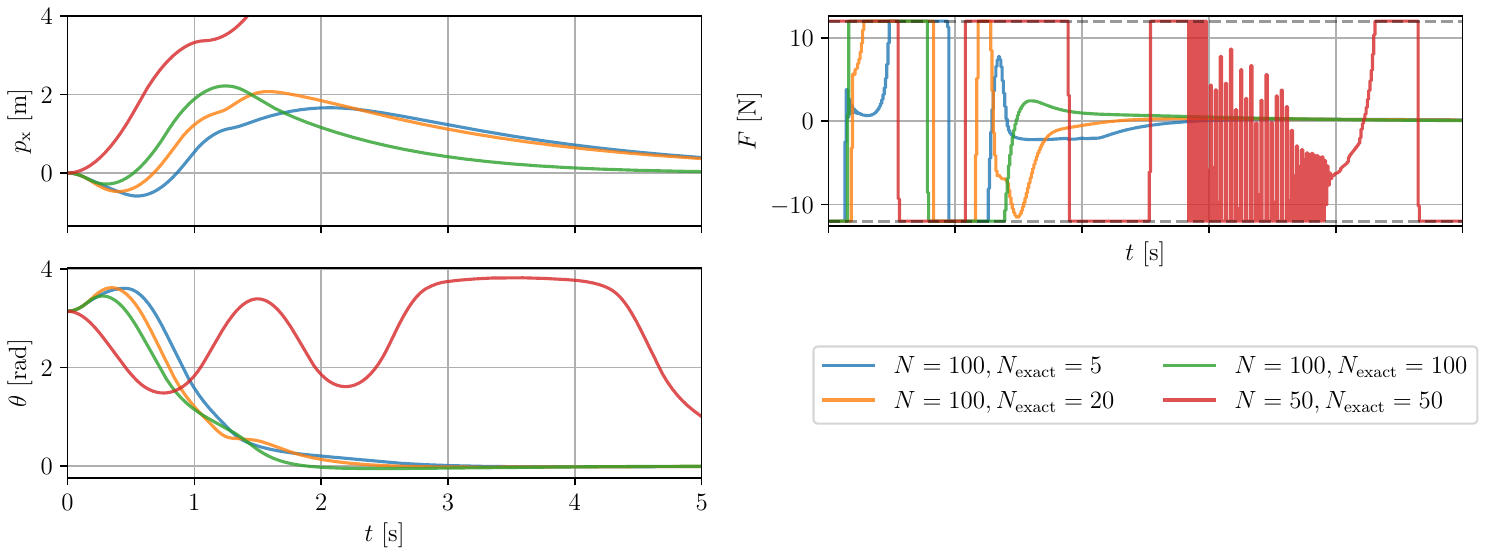}
	\caption{Closed-loop trajectories corresponding to different controller variants with partial tightening. The total number of shooting intervals $N$ and the ones handled exactly $N\ind{exact}$ are varied. The shooting intervals are part of the tightened horizon.
		\label{fig:zanelli_ptight_trajectories}
	}
\end{figure}

This section briefly demonstrates how partial tightening can be efficiently implemented within \acados, qualitatively reproducing the results from the benchmark in the paper that introduced partial tightening~\cite{Zanelli2017a}.
To this end, we regard the inverted pendulum on cart model as in Section~\ref{sec:experiments_pendulum_mocp}.
The cost function is of linear least-squares form
\begin{align}
l(x, \contcontrl) = x\transp Q x + \contcontrl\transp R \contcontrl,
\end{align} 
with $Q = \diag{0.1, 1, 0.1, 2\cdot 10^{-3}}, R=5\cdot 10^{-4}$.
For the terminal cost, we use $ x\transp P x$, where $P$ is the solution of the continuous-time algebraic Riccati equation with cost and dynamics linearized at the steady-state.
The only constraint is that the control input is $-12 \leq \contcontrl \leq 12$.
When using partial tightening, the constraint is replaced with an additional log-barrier term as in~\eqref{eq:barrier} corresponding to this constraint with $\tau=5$ and a GGN Hessian approximation is used.
While previous works only used custom implementations of partial tightening, the MOCP interface allows a convenient formulation in established software and the source code to reproduce the results presented in this Section is publicly available\footnote{\url{https://github.com/FreyJo/ocp_solver_benchmark/blob/main/experiments/zanelli_partial_tightening.py}}.

This benchmark compares controllers in a closed-loop simulation with initial state $x_0 = (0, \pi, 0, 0)$ for a duration of $5\mathrm{s}$ with a time step of $0.01 \mathrm{s}$.
The controllers use different underlying OCP formulations, where the length of the overall horizon and the tightened horizon are varied.
The tightened horizon is implemented as a second phase of an MOCP formulation.
The length of a shooting interval is fixed to $0.01\mathrm{s}$ and the dynamics are discretized using an RK4 integrator.
We vary the total number of shooting intervals $N$ and the number of shooting intervals on which the constraints are formulated exactly $N\ind{exact}$.
The second phase contains $N-N\ind{exact}$ shooting intervals, on which the control bounds are replaced with logarithmic barriers.

All controllers use the real-time iteration algorithm.
The QPs are solved using partial condensing, such that the blocks corresponding to the tightened horizon are condensed into one block and the remaining blocks remain uncondensed.
In order to perform as many operations as possible in the preparation phase, \acados~internally implements functions that assume that only matrices of the QP are known and a second one that completes the computations once the vector quantities are known.
For the \acados{} module that performs the condensing and QP solution has a split functionality, namely \texttt{condense\_lhs()} and \texttt{condense\_rhs\_and\_solve()}.

Table~\ref{tab:ptight_performance} gives an overview on the controllers closed-loop performance and the computation times split into preparation and feedback phase.
Figure~\ref{fig:zanelli_ptight_trajectories} visualizes the closed-loop trajectories for different controller variants.
Figure~\ref{fig:zanelli_ptight_rti_timings} and Table~\ref{tab:ptight_performance} show that the preparation time is consistent between all solver variants with $N=100$ and roughly halved for the variant with $N=50$.
While the controller with $N=50, N\ind{exact} = 50$ fails at swinging up the pendulum, the other controllers succeed at this task.
Comparing the variants with $N=100$, one can see that decreasing $N\ind{exact}$ results in an increase in relative suboptimality and a decrease in associated computational complexity, more specifically a decrease int he computation time of the feedback phase.
Overall, the results are very similar to the ones reported in the original benchmark\cite{Zanelli2017a}.
Note that in our implementation, the barrier parameter $\tau$ could be easily increased over the shooting intervals, resulting in a progressive tightening formulation that is not just partial tightening.
In summary, partial and progressive tightening formulations allow one to trade-off computational complexity and closed-loop performance and can be formulated conveniently as MOCPs.

\FloatBarrier
\section{Conclusion}
This paper gives an overview on multi-phase OCP (MOCP) formulations and their efficient treatment using multiple shooting.
Several approaches are presented which allow to formulate a continuous-time OCP in a successively approximate way which require an MOCP formulation.
The work provides an overview on different control parametrizations for use within multiple shooting, such as piecewise polynomial controls of different degrees and closed-loop costing variants, and motivated their use.
These control parametrizations have been compared on a benchmark example from previous work.
Moreover, we demonstrate the efficiency of NMPC controllers based on an MOCP formulation with an approximate model in a second phase.
Lastly, we show how partial and progressive tightening OCPs can be phrased as MOCPs.
These examples show that the added degrees of freedom allow one to develop NMPC controllers, which outperform ones limited to single-phase OCP formulations.
While other state-of-the-art software packages are limited to single-phase OCP formulations or use general purpose NLP solvers instead of structure-exploiting algorithms tailored to OCPs, the new \acados~feature allows for both a convenient formulation and the generation of efficient solver for MOCPs.
We believe that this new feature will make efficient solvers for MOCP formulations more available to NMPC practitioners and spread their use in real-world applications.

\bibliographystyle{ieeetran}
\bibliography{syscop}

\end{document}